\documentclass[12pt]{amsart}
\usepackage{amsmath, amsthm, amssymb,cite,enumitem}
\usepackage{fullpage}
\usepackage{color}
\usepackage{hyperref}
\usepackage[english]{babel}
\usepackage{framed}
\usepackage{enumitem}

\newtheorem{theorem}{Theorem}
\newtheorem{lemma}{Lemma}
\newtheorem{proposition}[lemma]{Proposition}
\newtheorem{corollary}[lemma]{Corollary}

\newtheorem{remark}[lemma]{Remark}
\newtheorem{conjecture}[lemma]{Conjecture}
\numberwithin{lemma}{section}

\numberwithin{equation}{section}

\newcommand{\R}{{\mathbb R}}

\newcommand\Tau{\mathcal{T}}

\usepackage[makeroom]{cancel}

\newcommand{\sgn}{\mathop{\mathrm{sgn}}}

\renewcommand\R{\ R}

\newcommand{\cB}{\mathring{B}}
\newcommand{\cH}{\mathring{H}}

\begin{document}

\title{The lifespan of small data solutions for Intermediate Long Wave equation (ILW)}

\author{Mihaela Ifrim}
\address{Department of Mathematics, University of California at Berkeley}
\thanks{The first author was supported by the Sloan Foundation, and by an NSF CAREER grant DMS-1845037. }
\email{ifrim@math.wisc.edu}

\author{Jean-Claude Saut}
\address{Laboratoire de Mathématiques \& Université de Paris - Saclay, 91405 Orsay, France}
\email{jean-claude.saut@universite-paris-saclay.fr}

\begin{abstract}
  This article represents a first step toward understanding the long
  time dynamics of solutions for the Intermediate Long Wave equation (ILW).  While
  this problem is known to be both completely integrable and globally
  well-posed in $H^{\frac{3}{2}}$, much less seems to be known concerning its long
  time dynamics. Here we prove well-posedness at much lower regularity, namely an $L^2$ global well-posedness result. Then we consider the case of  small and localized  data and show that the solutions disperse up to cubic timescale.
\end{abstract}

\maketitle

\section{Introduction}

In this article we consider  the  Intermediate Long Wave equation (ILW)
\begin{equation}\label{ilw}
\left\{
\begin{aligned}
&(\partial_t +\frac{1}{\delta}\partial_x+\Tau_{\delta }^{-1} \partial_x^2) \phi  =   \frac12 \partial_x (\phi^2)\\
&\phi(0) = \phi_0,\\
\end{aligned}
\right.
\end{equation}
where $\phi$ is a real valued function, $\phi :\mathbb{R}\times \mathbb{R}\rightarrow \mathbb{R}$, and  $\Tau_\delta$ is the Tilbert transform with a large scale parameter $\delta$. Precisely $\Tau^{-1}_{\delta}$ is a zero order operator with symbol $i \coth (\delta \xi)$, arising in the description of the Dirichlet to Neumann map in a two dimensional finite depth domain occupied by a two-layer  fluid, see for instance \cite{BLS}.
In the limit as $\delta \rightarrow \infty $ this converges to the Hilbert transform, which corresponds to the  Dirichlet to Neumann map in  infinite depth.

The ILW equation appears as a model for long internal gravity waves in stratified fluids, and $\delta$ plays the role of a relative depth 
factor, roughly measuring the ratio of the deeper, lower fluid depth  and the top, shallower fluid depth. As $\delta$ approaches infinity,
one obtains the Benjamin-Ono equation in the limit.

 The Hilbert transform is given by 
\begin{equation}
\label{hilbert}
\mathcal{H}\phi(x)=-\frac{1}{\pi}\, p.v. \int_{-\infty}^{\infty} \frac{\phi(y)}{x-y}\, dy,
\end{equation}
whereas  $\Tau^{-1}_{\delta}$ is defined by the principal value convolution
\begin{equation}
\label{till}
\Tau^{-1}_{\delta}\phi(x)=\frac{1}{2\delta}\, p.v. \int_{-\infty}^{\infty} \coth \left( \frac{\pi (x-y)}{2\delta}\right)\phi(y)\, dy.
\end{equation}

This can be seen as an inverse of the Tilbert transform which is defined as follows:
\[
\Tau_{\delta} f(x)=-\frac{1}{2\delta} p.v. \int_{-\infty}^{\infty} \mbox{cosech} \left( \frac{\pi}{2\delta} (x-x')\right) f(x')\, dx'.
\]
The convention here is that the intermediate long wave equation models  unidirectional waves that travel to the right; which in turn determines the sign choices in the equation \eqref{ilw}. This corresponds to 
a nonnegative group velocity for the associated linear flow.

The operator  $\Tau_{\delta}$ vanishes on constant functions, therefore there is an ambiguity in the definition of $\Tau^{-1}_{\delta}$, which is apriori defined only modulo constants but this is not an issue as it always comes paired with at least one spatial derivative $\partial_x$. For concreteness we define it to be given by the Fourier multiplier p.v.\,$ i\coth (\delta\xi)$, which is consistent with the kernel given above.

The intermediate long wave (ILW) equation was first introduced in \cite{j} and \cite{kkd}. There, the reader will find a  more detailed history of how the equation was derived, as well as a comprehensive discussion on the  physical relevance of this one dimensional wave model equation. We refer to \cite{CGK, BLS} for a rigorous derivation (in the sense of consistency) from the two-layer system for internal waves, see also \cite{Xu} for a rigorous derivation of the system version of the ILW equation in one or two spatial variables,  including the case of a free upper surface. 

 The ILW equation has received a fair amount of attention lately, and there are multiple theoretical works devoted to its analysis. We refer to the survey paper \cite{klein-saut} and to  the book \cite{KS2} Chapter 3 for an exhaustive description of the relevant works.

We  observe that one can make the change of variable $x\rightarrow x-\delta^{-1} t$ and transform \eqref{ilw} into
\begin{equation}\label{ilw-}
\left\{
\begin{aligned}
&(\partial_t +\Tau_{\delta }^{-1} \partial_x^2) \phi  =   \frac12 \partial_x (\phi^2)\\
&\phi(0) = \phi_0.\\
\end{aligned}
\right.
\end{equation}
Thus, the linear transport term $\delta^{-1} \phi_x$ is artificial. It is a consequence of the different scalings arising in the regimes where the ILW equations appears: in the surface layer, the long-wave scaling needs to be matched to a different scaling in the deep
lower layer, where the vertical scale matches the horizontal scale. We also observe that  the term  $\delta^{-1} \phi_x$  insures that  the model is centered around waves with zero velocity (minimum speed), which emphasizes the connection with the  Korteweg–De Vries (KdV) equation in the low frequency regime.  
\begin{remark}\label{r:transport} In this paper we will use this change of coordinates in the proof of the $L^2$ global well-posedness result in order to emphasize the connection with the Benjamin-Ono equation. On the other hand, we will keep the transport term in the proof of the dispersive result in order to emphasize the connection with the KdV equation.

\end{remark}

One final remark about $\delta$ is that it can be thought of as a scaling parameter. Indeed, one may set $\delta =1$ with a linear change of coordinates and substitution
\[
v(t,x) = \delta u(\delta^2 t,\delta x).
\]
From here on, we will work with $\delta = 1$ in the present paper. 
\bigskip

Equation \eqref{ilw} is known to be completely integrable, \cite {KSA, KAS}. In particular it possesses a Lax pair and   an infinite hierarchy of
conservation laws. However a rigorous theory of the Cauchy problem using Inverse Scattering techniques is not available so far. Recent progress has been made for the direct scattering problem with small data in weighted $L^2$ spaces, see \cite{Klip, KPW}. Anyway integrability  cannot be at the moment connected with the Cauchy
problem at low regularity.

We list only some of these conserved energies, which hold for smooth solutions (for example $H_x^3(\mathbb{R})$).  
Integrating by parts, one sees that this problem has conserved mass,
\[
E_0 = \int \frac{1}{2} \phi^2 \, dx,
\]
momentum
\[
E_1 = \int \phi \Tau^{-1}\phi_x - \frac13 \phi^3\, dx ,
\]
as well as energy
\[
E_2 = \int \frac{1}{2}\phi_x^2 -\frac32 \phi^2 \Tau^{-1} \phi_x + \frac14 \phi^4 +\frac{3}{2}\left[\Tau^{-1}\phi_x
\right]^2\, dx.
\]

Since the ILW equation is completely integrable, at each nonnegative integer $k$
we similarly have a conserved energy $E_k$, see for instance \cite{abfs} for the explicit form of $E_3.$ In order to be able to verify that these energies are conserved, (including the ones not listed here), one needs to use the well known identity of the integral operator $\Tau$:
\begin{equation}
\label{tau ident}
\Tau^{-1} (u\Tau v + v\Tau u)=uv-(\Tau u)(\Tau  v).
\end{equation}

Considering local and global well- posedness results in Sobolev spaces $H^s$, a natural threshold is given by the scale invariance law that the ILW equation satisfies in the high frequency limit:
\begin{equation}
\label{scaling}
\phi(t,x)\rightarrow \lambda \phi (\lambda^2t,\lambda x),
\end{equation}
and the scale invariant Sobolev space associated to this scaling is $\dot{H}^{-\frac{1}{2}}$. This is the exact scale invariant Sobolev space we have for the Benjamin-Ono equation.

There have not been many developments in the well-posedness theory of the Intermediate Long Wave equations. Its close connection with the study of the two dimensional water waves equations makes it a very interesting model to understand. An extensive discussion of the ILW equation and related models can be found in the survey papers \cite{klein-saut}, and in the book \cite{KS2}.

We are concerned with the Cauchy problem at low regularity for our equation \eqref{ilw}. The first global well-posedness  result at the level of $H^s(\mathbb{R})$ was obtained in \cite{abfs} for the Sobolev index $s> 3/2$. The best known results are due to Molinet and Vento , \cite{MV} for $s\geq 1/2,$ the unconditional uniqueness holding for $s>1/2,$ respectively to Molinet-Pilod-Vento \cite{Molinet2017OnWF}  where well-posedness for ILW is achieved in $H^s$, $s>\frac{1}{4}$. Well-posedness in the range $-\frac12 \leq s< \frac{1}{4}$ appears to be an open question. 

Very few results are known concerning the global behavior of solutions. In \cite{MuPoSa}, extending previous results in \cite{MuPo} for the Benjamin-Ono equation, the authors establish local energy decay  in an increasing in time region $I(t)$ of space of size $t/\log (t )$ implying in particular the nonexistence of breathers inside the region $I(t)$ for any time t sufficiently large.

\vspace{0.3cm}
Our goal in the present paper is two fold: 

\begin{enumerate}[label=\roman*)]

\item We extend the  well-posedness range of equation \eqref{ilw} to  Sobolev indices  $0\leq s\leq \frac{1}{4}$ , and 

\item We prove dispersive decay estimates  for the nonlinear problem  on the cubic time scale. 
\end{enumerate}

We begin our discussion with the local well-posedness problem, where we first review some key thresholds in the analysis.
The $H^s$ with $s\geq 3/2$ well-posedness result was obtained in
\cite{abfs} (with uniqueness for $s>3/2$) using energy estimates and an interpolation argument. For convenience we use this result as a starting point for
our work, which is why we recall it here:

\begin{theorem} 
The ILW equation is globally  well-posed in $H^s$, $s\geq 3/2$.
\end{theorem}

Our first goal here is to provide an $L^2$ theory for the ILW equation,
and prove the following theorem:
\begin{theorem}\label{thm:lwp}  
The ILW equation is globally well-posed in $L^2$.
\end{theorem}

Since the $L^2$ norm of the solutions is conserved, this is in effect a local in time result,
trivially propagated in time by the conservation of mass. In particular, it says little 
about the long time properties of the flow, which will be our primary target here.

Given the quasilinear nature of the ILW equation, here it is important 
to specify the meaning of well-posedness. This is summarized in the following properties:

\begin{description}
\item [(i) Existence of regular solutions] For each initial data $\phi_0 \in H^{\frac32}$ there exists a unique 
global solution $\phi \in C(\mathbb{R};H^{\frac32})$.

\item[(ii) Existence and uniqueness of rough solutions] For each initial data $\phi_0 \in L^2$ 
there exists a solution $\phi \in  C(\mathbb{R};L^2)$, which is the unique limit of regular solutions.

\item[(iii) Continuous dependence ] The data to solution map $\phi_0 \to \phi$ is continuous 
from $L^2$ into $C(L^2)$, locally in time.

\item [(iv) Higher regularity] The data to solution map $\phi_0 \to \phi$ is  continuous 
from $H^s$ into $C(H^s)$, locally in time, for each $s > 0$.

\item[(v) Weak Lipschitz dependence] The flow map for $L^2$ solutions is locally Lipschitz
in the $ L^2$ topology. 
\end{description}

We recall, see \cite{MST} that the Intermediate Long Wave equation is {\it quasilinear} in the sense that the flow map cannot be smooth, say $C^2.$

The ideas we will use to prove the above theorem are directly inspired by the work of Ifrim-Tataru \cite{IT-BO} where, among other results, low regularity global well-posedness result and weak Lipschitz dependence for the Benjamin-Ono equation were obtained.

In a nutshell, the  approach to this result is based on the idea of normal forms, introduced by Shatah \cite{shatah} in the dispersive realm in the context of studying the long time behavior of dispersive pde's. Here we turn it around and consider it in the context of studying local well-posedness.  In doing this, the main difficulty we face is that the standard normal form method does not readily apply for quasilinear equations. To resolve this difficulty, we will use some ideas which were developed in the Benjamin-Ono context in \cite{IT-BO}.

\smallskip

Our second goal is to establish dispersive bounds for the solution in the case of  small  and localized initial data.  Given a
generic quasilinear problem with data of size $\epsilon$ and quadratic interactions, the standard result is to obtain quadratic lifespan bounds, i.e., $T_{max} \lesssim \epsilon^{-1}$. This does not require any localization. But if the initial data is also localized, then one may hope to show that the solutions also exhibit dispersive decay, possibly on a better time scale.

Here we show that for our problem, despite the presence of quadratic interactions, the time threshold for dispersive decay  is nevertheless at least cubic, i.e., $T_{max} \gtrsim \epsilon^{-2}$.

\begin{theorem}
\label{t:quartic} Consider the ILW equation \eqref{ilw-} with small and localized initial data $\phi_{0}$,
\[
\Vert \phi_0\Vert_{L^2} +\Vert x\phi_0\Vert_{L^2}\lesssim \epsilon. 
\]
Then the solution $\phi$ exists and satisfies the dispersive decay bounds 
\begin{equation} \label{point1}
|\phi(t,x)| \lesssim
\left\{
\begin{aligned}
 & \epsilon t^{-\frac13} \left< t^{-1/3}x \right>^{-1/4} \left<t^{-1/3}x_{+} \right>^{-3/4+} \mbox{ for } x \geq - t \\
 & \epsilon t^{-\frac{1}{2}} \mbox{ for } x < -t,
  \end{aligned}
  \right.
\end{equation}
and
\begin{equation} \label{point2}
\vert \Tau^{-1}\phi (t,x)\vert  \lesssim
\left\{
\begin{aligned}
 & \epsilon t^{-\frac23} \left< t^{-1/3}x\right>^{1/4} \left<t^{-1/3}x_{+} \right>^{-5/4} \mbox{ for } x \geq  - t\\
 & \epsilon t^{-\frac{1}{2}} \mbox{ for } x < - t .
  \end{aligned}
\right. 
\end{equation}
on a time interval $\left[ -T_{\epsilon}, T_{\epsilon} \right]$ with $T_{\epsilon}\lesssim \epsilon^{-2}$. 
\end{theorem}

To place this result into context, we recall that at high frequencies the ILW evolution closely resembles the Benjamin-Ono
equation 
\begin{equation}\label{BO}
\left\{
\begin{aligned}
&(\partial_t +H \partial_x^2) \phi  =   \frac12 \partial_x (\phi^2)\\
&\phi(0) = \phi_0, \\
\end{aligned}
\right.
\end{equation}
while at low frequencies it is well approximated by the KdV flow 
\begin{equation}\label{IT-BO}
\left\{
\begin{aligned}
&(\partial_t + \partial_x^3) \phi  =   \frac12 \partial_x (\phi^2)\\
&\phi(0) = \phi_0. \\
\end{aligned}
\right.
\end{equation}

 To understand these comparisons it suffices to look at the dispersion relation of the linear counterpart of \eqref{ilw}
 (with $\delta = 1$),
\[
\tau =\xi ^2 \coth \xi -  \xi,
\]
which shows that  the ILW equation is a dispersive equation which at low frequencies resembles the  KdV-like dispersion
\[
\tau \approx \frac{1}{3} \xi^3,
\]
 and at high frequencies acts as a Benjamin Ono dispersion
 \[
 \tau \approx \xi|\xi| -\xi,
 \]
where the transport term can be eliminated via a linear change of coordinates.
 
The heuristics here are that in some sense ILW equation borrows features from both equations, and thus we encounter difficulties specific both to NLS and to KdV respectively. 

The dispersive decay of solutions for the Benjamin-Ono equation
was considered in \cite{IT-BO}, where it was shown that 
the small data solutions have dispersive decay up to an almost global time $T_{BO} \approx e^{-\frac{c}{\epsilon}}$.  On the other
hand the dispersive decay of solutions for the KdV equation
was considered in \cite{KIT}, where it was shown that 
the small data solutions have dispersive decay up to a quartic time $T_{KdV} \approx \epsilon^{-3}$. Both of these results are optimal, as the above time scales can be identified heuristically with the earliest possible emergence time of a solitary wave arising from the data; this uses inverse scattering, and is based on a spectral analysis of the corresponding Lax operator. 

One also has solitary wave solutions for the ILW equation,
which resemble the KdV ones at low amplitude (which are the ones that may emerge from small data) and the Benjamin-Ono 
ones at high amplitude. An important role here is played by the 
symbol smoothness: the ILW symbol is smoother than the Benjamin-Ono symbol, and  such a property is responsible for  many qualitative differences in the behaviour of the solutions of \eqref{ilw}: e.g. explicit  N-soliton  solutions decay exponentially at infinity, contrary to those of the Benjamin-Ono equation, see for example \cite{m, j}.

Given the above discussion, we would expect that the 
optimal time scale for the ILW to also be quartic:

\begin{conjecture}
The result of theorem~\ref{t:quartic} should hold up to a quartic time $t \lesssim \epsilon^{-3}$.
\end{conjecture}

However, proving such a result in the ILW case is considerably more complex than the KdV case, and would require improvements at every single step
of the way.

The key element of the proof of Theorem~\ref{t:quartic} is to 
obtain a good (nonlinear) vector field bound. In the Benjamin-Ono case, it was discovered in \cite{IT-BO} that there exists in effect a vector field type conserved quantity, which led to proving almost global in time dispersive properties for the solutions. This is no longer the case for KdV. There a good vector field is instead suggested
by the scaling symmetry, but it only yields an approximate 
conservation law. Neither of these ideas applies directly in the 
ILW case, so here we need to construct a good vector field. Indeed,
most of the work in Section 5 is devoted to this construction.
We also refer the reader to \cite{IT-wp}, where a related construction was done for problems with cubic nonlinearities.

The set of resonances for ILW is  important in our vector field construction.  The resonant two wave interactions correspond to solutions to the system
 \[
 \left\{
 \begin{aligned}
 &(\zeta)^2\coth \zeta= \xi^2 \coth \xi +\eta^2\coth \eta\\
 &\xi +\eta=\zeta.
  \end{aligned}
  \right.
 \]
 The only solutions occur when at least one of $\xi, \eta, \zeta$ vanishes. Formally,  quadratic resonances can be removed by means of a normal form transformation of the form
\[
\tilde{\phi}=\phi +\mathcal{B}(\phi, \phi),
\]
where the symbol of $\mathcal{B}$ is 
\[
\mathcal{B}(\xi, \eta) = -\frac{i}{2}\frac{(\xi+\eta)}{\left[ (\xi +\eta)^2 \mbox{coth}(\xi +\eta) -\eta^2 \mbox{coth}(\eta) -\xi^2\mbox{coth}(\xi)\right]} \, .
\]
This is nonsingular outside of the resonance set, but has singularities on the resonance set. The key to our result is that these singularities no longer appear in the vector field construction.
This is  closely related to a robust adaptation of the normal form method to quasilinear equations, called \emph{quasilinear modified energy method}. It was
introduced earlier by Ifrim-Hunter-Tataru-Wong in \cite{BH}, and
then further developed in the water wave context first in \cite{HIT}
and later in \cite{IT-c, itg, IT-g, itc, IT-BO}. There the idea is to modify the energies, rather
than apply a normal form transform to the equations; this method  is then
successfully used in the study of long time behavior of solutions.
Alazard and Delort~\cite{alazarddelort, ad1} have also developed another way of constructing 
the same type of almost conserved energies by using a partial normal form transformation
to symmetrize the equation, effectively diagonalizing the leading part of the energy. 

\subsection*{Acknowledgements:} The first author was supported by the Sloan Foundation, and by  NSF CAREER grant DMS-1845037.

\section{Definitions and review of notations }

\subsubsection*{The big O notation:} We use the notation $A \lesssim
B$ or $A = O(B)$ to denote the estimate $|A| \leq C B$, where $C$
is a universal constant which will not depend on $\epsilon$.  If $X$
is a Banach space, we use $O_X(B)$ to denote any element in $X$ with
norm $O(B)$; explicitly we say $u =O_X(B)$ if $\Vert u\Vert _X\leq C
B$.  We use $\langle x\rangle$ to denote the quantity $\langle x
\rangle := (1 + |x|^2)^{1/2}$.

\subsubsection*{Littlewood-Paley decomposition:} One important tool in
dealing with dispersive equations is the Littlewood-Paley
decomposition.  We recall its definition and also its usefulness in
the next paragraph. We begin with the Riesz decomposition
\[
 1 = P_- + P_+,
\]
where $P_\pm$ are the Fourier projections to $\pm [0,\infty)$; from 
\[
\widehat{Hf}(\xi)=-i\sgn (\xi)\, \hat{f}(\xi),
\]
 we observe that
\begin{equation}\label{hilbert-re}
iH = P_{+} - P_{-}.
\end{equation}
Let $\psi$ be a bump function adapted to $[-2,2]$ and equal to $1$ on
$[-1,1]$.  We define the Littlewood-Paley operators $P_{k}$ and
$P_{\leq k} = P_{<k+1}$ for $k \geq 0$ by 
$$ \widehat{P_{\leq k} f}(\xi) := \psi(\xi/2^k) \hat f(\xi)$$
for all $k \geq 0$, and $P_k := P_{\leq k} - P_{\leq k-1}$ (with the
convention $P_{\leq -1} = 0$).  Note that all the operators $P_k$,
$P_{\leq k}$ are bounded on all translation-invariant Banach spaces,
thanks to Minkowski's inequality.  We define $P_{>k} := P_{\geq k-1}
:= 1 - P_{\leq k}$.

For simplicity, and because $P_{\pm}$ commutes with the
Littlewood-Paley projections $P_k$ and $P_{<k}$, we will introduce the
following notation $P^{\pm}_k:=P_kP_{\pm}$ , respectively
$P^{\pm}_{<k}:=P_{\pm}P_{<k}$.  In the same spirit, we introduce the
notations $\phi^{+}_k:=P^{+}_k\phi$, and $\phi^{-}_k:=P^{-}_k\phi$,
respectively.

Given the projectors $P_k$, we also introduce additional projectors $\tilde P_k$
with slightly enlarged support (say by $2^{k-4}$) and symbol equal to $1$ in the support of $P_k$.

From Plancherel's theorem  we have the bound
\begin{equation}\label{eq:planch}
 \| f \|_{H^s_x} \approx (\sum_{k=0}^\infty \| P_k f\|_{H^s_x}^2)^{1/2}
\approx (\sum_{k=0}^\infty 2^{2ks} \| P_k f\|_{L^2_x}^2)^{1/2}
\end{equation}
for any $s \in \mathbb{ R}$.

\subsubsection*{Multi-linear expressions} We shall now make use of a
convenient notation for describing multi-linear expressions of product
type, as in \cite{tao-WM}.  By $L(\phi_1,\cdots,\phi_n)$ we denote a
translation invariant expression of the form
\[
L(\phi_1,\cdots,\phi_n)(x) = \int K(y) \phi_1(x+y_1) \cdots \phi_n(x+y_n) \, dy, 
\] 
 where $K \in L^1$. More generally, one can replace $Kdy$ by any bounded measure. 
  By $L_k$ we denote such multilinear expressions whose output is localized at frequency $2^k$.
 
  This $L$ notation is extremely handy for expressions such as the
  ones we encounter here; for example we can re-express the normal
  form \eqref{commutator B_k} in a simpler way as shown in
  Section~\ref{s:local}. It also behaves well with respect to reiteration,
e.g. 
\[
L(L(u,v),w) = L(u,v,w).
\]

Multilinear $L$ type expressions can easily be estimated
in terms of linear bounds for their entries. For instance we have
\[
\Vert L(u_1,u_2) \|_{L^r} \lesssim \|u_1\|_{L^{p_1}} \|u_2\|_{L^{p_2}}, \qquad \frac{1}{p_1}+ \frac{1}{p_2} = \frac{1}{r}.
\]
A slightly more involved situation arises in this article when we seek to 
use bilinear bounds in estimates for an $L$ form. There we need to account for the 
effect of uncorrelated translations, which are allowed given the integral bound 
on the kernel of $L$. To account for that we use the translation group $\{T_y\}_{y \in \R}$,
\[
(T_y u)(x) = u(x+y),
\]
and estimate, say, a trilinear form as follows:
\[
\|L(u_1,u_2,u_3) \|_{L^r} \lesssim \|u_1\|_{L^{p_1}} \sup_{y \in \R} 
\|u_2 T_y u_3\|_{L^{p_2}}, \qquad \frac{1}{p_1}+ \frac{1}{p_2} = \frac{1}{r} .
\]
On occasion, we will write this in a shorter form
\[
\|L(u_1,u_2,u_3) \|_{L^r} \lesssim \|u_1\|_{L^{p_1}}   \|L(u_2,u_3)\|_{L^{p_2}}.
\]

  To prove the boundedness in $L^2$ of the normal form transformation,
  we will use the following proposition from Tao \cite{tao-WM}; for
  completeness we recall it below:

\begin{lemma}[Leibnitz rule for $P_k$]\label{commutator}  We have the commutator identity
\begin{equation}
\label{commute}
\left[ P_k\, ,\, f\right]  g = L(\partial_x f, 2^{-k} g).
\end{equation}
\end{lemma}
When classifying cubic terms (and not only) obtained after
implementing a normal form transformation, we observe that having a
commutator structure is a desired feature.  In particular Lemma~\ref{commutator} tells us that when one of the entry (call it $g$) has
frequency $\sim 2^k$ and the other entry (call it $f$) has frequency
$\lesssim 2^k$, then $P_k(fg) - f P_k g$ effectively shifts a
derivative from the high-frequency function $g$ to the low-frequency
function $f$.  This shift will generally ensure that all such
commutator terms will be easily estimated.

\begin{lemma}[Properties of the $ \mathcal{P}:=H-\Tau^{-1}$ operator]\label{l:hilbert-tilbert} The zero order operator $ \mathcal{P}:=H-\Tau^{-1}$ is a smoothing operator that acts like an antiderivative close to the zero frequency and decays very fast at high frequencies 
\begin{equation}
\label{nice}
\Vert  \mathcal{P} \partial^n_xf \Vert_{L^p} \lesssim \Vert f\Vert _{L^p}, \quad \mbox{ for any } p\geq 2 \mbox{ and } n\geq 1. 
\end{equation}
\end{lemma}
\begin{proof}
It suffices to inspect the symbol of the operator $\mathcal{P}$; its boundedness and smoothness are enough to conclude the $L^p\rightarrow L^{p}$, $p\in [2, \infty)$ mapping property. We have that the symbol of $\mathcal{P}$ is given by
\[
p(\xi)= -i (\sgn \xi  + \coth \xi ) = -i \left[ \pm 1 \pm  \frac{2}{e^{2\xi} -1}\right].
\]
At $\xi =0$, the symbol $p$ has a singularity but $\mathcal{P}$ is always paired with at least one derivative. In our current work,  the right object that will be on interest is $\mathcal{P}\partial^2_x$ which has the symbol $-\xi^2 p(\xi)$, which is a $C^{1,1} (\mathbb{R})$ object that go to zero exponentially. Same observation holds true for $\mathcal{P}\partial^n_x$ for any $n\geq 1$, meaning the symbol is in $C^{n-1,1}(\mathbb{R})$ and exponentially decaying. Thus, one gets the desired bound \eqref{nice}.

\end{proof}
 
\subsubsection*{Frequency envelopes.} In preparation for  the proof of one of  the main
theorems of this paper, we revisit the \emph{frequency envelope}
notion; it will turn out to be very useful, and also an elegant tool
used later in the proof of the local well-posedness result. Precisely, it appears in both (i) the
proof of the a-priori bounds for solutions for the Cauchy problem
\eqref{BO} with data in $L^2$, which we state in
Section~\ref{s:local}, and  (ii) the proof of the bounds for the
linearized equation, in Subsection~\ref{s:linearized}.

Following Tao's paper \cite{tao}, we say that a sequence of nonnegative real $c_{k}\in l^2$
is an $L^2$ frequency envelope for $\phi \in L^2$ if
\begin{itemize}
\item[i)] $\sum_{k=0}^{\infty}c_k^2 \lesssim 1$;\\
\item[ii)] it is slowly varying, $c_j /c_k \leq 2^{\delta \vert j-k\vert}$, with $\delta$ a universal constant;\\
\item[iii)] it bounds the dyadic norms of $\phi$, namely $\Vert P_{k}\phi \Vert_{L^2} \leq c_k$. 
\end{itemize}
Given a frequency envelope $c_k$ we define 
\[
 c_{\leq k} = (\sum_{j \leq k} c_j^2)^\frac12, \qquad  c_{\geq k} = (\sum_{j \geq k} c_j^2)^\frac12.
\]

\begin{remark}
To avoid dealing with certain issues arising at low frequencies,
we can harmlessly make the extra assumption that $c_{0}\approx 1$.
\end{remark}

\begin{remark}
Another useful variation is to weaken the slowly varying assumption to
\[
2^{- \delta \vert j-k\vert} \leq    c_j /c_k \leq 2^{C \vert j-k\vert}, \qquad j < k,
\]
where $C$ is a fixed but possibly large constant. All the results in this paper are compatible with this choice.
This offers the extra flexibility of providing higher regularity results by the same argument.
\end{remark}

\section{The linear flow}
Here we consider the linear ILW flow with $\delta = 1$,
\begin{equation}\label{ilw-lin}
(\partial_t + \partial_x+ \Tau^{-1} \partial^2_x)\psi = 0, \qquad \psi(0) = \psi_0,
\end{equation}
which we rewrite as 
\[
[ \partial_t -i A(D)]\psi =0,
\]
where 
\[
A(D)= i(\partial_x +\Tau^{-1}
\partial^2_x) 
\]
has real, odd symbol
\[
A(\xi) = \xi^2 \coth  \xi- \xi.
\]
The solution $\psi(t) = e^{it A(D)} \psi_0$ has conserved $L^2$ norm. 

Our aim in this section is to discuss the linear dispersive properties for this flow. We do this in two stages. First we consider uniform dispersive decay bounds 
under the assumption that the initial data is localized.
Then we consider $L^2$ type data, in which case the dispersion is best captured by Strichartz and bilinear $L^2$ bounds.

 \subsection{Uniform dispersive bounds.}
 To write these we define the weights
\begin{equation} \label{def-omega0}
 \omega_0(t,x) = 
 \left\{
 \begin{aligned}
 & t^{-\frac13} \left< t^{-1/3}x \right>^{-1/4} \left<t^{-1/3}x_{+} \right>^{-3/4+} \mbox{ for } x \geq - t \\
 &t^{-\frac{1}{2}}   \qquad  \qquad \qquad  \qquad \qquad  \qquad  \, \, \mbox{   for } x < -t,
  \end{aligned}
  \right.
\end{equation}
\begin{equation} \label{def-omega1}
 \omega_1(t,x) = 
 \left\{
 \begin{aligned}
 &t^{-\frac23} \left< t^{-1/3}x\right>^{1/4} \left<t^{-1/3}x_{+} \right>^{-5/4} \, \mbox{ for } x \geq  - t\\
 &t^{-\frac{1}{2}} \ \qquad  \qquad \qquad  \qquad \qquad   \quad \mbox{ for } x < - t  .
  \end{aligned}
  \right.
\end{equation}
Here $x_+$ denotes the positive part of $x$, and the extra decay when $x > 0$ corresponds to the fact that linear ILW  waves travel to the left, which is similar to the standard KdV model. We begin our discussion with a standard 
dispersive decay bound, namely

\begin{proposition}\label{p:bodisp}
The fundamental solution 
\[
K(t,x) = e^{-tA} \delta_0 
\]
to the linear ILW flow satisfies the dispersive bounds
\begin{equation}
\begin{aligned}
 |K(t,x)| \lesssim & \  \omega_0(t,x)  , 
\\
|\Tau(D) K(t,x)| \lesssim & \ \omega_1(t,x).
\end{aligned}
\end{equation}
\end{proposition}

This is a well known result, which is proved using the stationary phase method, since the kernel $K$ 
admits the Fourier representation
\[
K(t) = \mathcal F^{-1} e^{ita(\xi)}.
\]
One may improve the above bounds to exponential decay in the right quadrant, but this does not play 
a significant role on what follows. The above decay rate depends on the velocity, at least at low velocity 
which corresponds to low frequency. Because of this, it is also interesting to 
consider a frequency localized version of the above proposition. Here we only consider localizations to low dyadic frequencies $2^j$ with $j \leq 0$, decomposing the fundamental solution $K$  as 
\[
K = P_{> 0} K + \sum_{j \leq 0} P_j K : = K_{>0} + \sum_{j < 0}  K_j.
\]

\begin{proposition}\label{p:bodisp-loc}
The dyadic portions $K_j$, $K_{> 0}$ of the fundamental solution $K$  
for the linear ILW flow satisfy the dispersive bounds
\begin{equation}\label{K0}
\begin{aligned}
 |K_{> 0}(t,x)| + |\Tau(D) K_{>0}(t,x)| \lesssim 
\left\{
 \begin{aligned} 
 & t^{-\frac12}  \qquad \mbox{ for } x < - \frac14 t \\
 &t^{-\frac{1}{2}} (1+|x|+t)^{-N} \mbox{ otherwise, }
  \end{aligned}
 \right.
\end{aligned}
\end{equation}
respectively
\begin{equation}\label{Kj}
\begin{aligned}
 |K_j(t,x)| + 2^{-j} |\Tau(D) K_{j}(t,x)| \lesssim 
\left\{
 \begin{aligned} 
 & 2^{-\frac{j}2} (t+2^{-3j})^{-\frac12}  \qquad \mbox{ for }  -  2^{2j+2} t <  x < -  2^{2j-2} t \\
 & 2^{-\frac{j}2} (t+2^{-3j})^{-\frac12} (1+2^j|x|+2^{3j} t)^{-N} \mbox{ otherwise. } 
  \end{aligned}
 \right.
\end{aligned}
\end{equation}
\end{proposition}

This is also easily proved using stationary phase, and reflects the fact that waves with frequency $> 1$ 
move with speed $\lesssim -1$, and waves with frequency $2^j$ move with speed $\approx -2^{2j}$.
We remark here that Proposition~\ref{p:bodisp} can be directly obtained from Proposition~\ref{p:bodisp-loc}
simply by dyadic summation. 

Our main interest in this paper is in the nonlinear ILW flow, where we cannot use Fourier methods anymore. 
So instead, we will reinterpret these decay bounds from an $L^2$ perspective.
We will denote by $P$ the linear operator associated to the ILW equation,
\[
P:= \partial_t -itA(D).
\]
To measure initial data localization we will use the operator $L$ given by
\[
L:=x+ tA_{\xi}(D),
\]
which is the push forward of $x$ along the linear flow,
\[
L(t) = e^{it A(D)} x e^{-it A(D)} ,
\]
and thus commutes with the linear operator,
\[
[ L, \partial_t - itA] = 0.
\]
Then for solutions $\psi$ to \eqref{ilw-lin},  the $L^2$ norm of $\psi$ 
and $L\psi$ are conserved. 

To make it useful in the nonlinear setting, we  will recast the dispersive bounds in Proposition~\ref{p:bodisp}
as a Sobolev type bound in terms of $u$ and $Lu$.  For this we need to make a good choice of Sobolev spaces, 
and this is inspired from prior work on the two limiting problems, namely KdV and Benjamin-Ono. In the Benjamin Ono case, see \cite{IT-BO}, it is most efficient to work with the $L^2$ norm for both $\psi$ and $L\psi$. 
On the other hand for KdV, it is better to use the space $B^{-\frac12}_{2,\infty}$ for $\psi$, respectively 
$\dot H^\frac12$ for $L\psi$. In our case, we will combine these two settings, using KdV style norms at low frequency and Benjamin-Ono style norms at high frequency. 
Hence we introduce Besov spaces $\cB^{-\frac12}_{2,\infty}$ respectively $\cH^{\frac12}$
with norms defined as follows:
\begin{equation}\label{B-12}
\| u \|_{\cB^{-\frac12}_{2,\infty}}^2 
:= \|u\|_{L^2}^2 + \sup_{k < 0} 2^{-k} \|P_k u\|_{L^2}^2,
\end{equation}
respectively 
\begin{equation}\label{H12}
\| v \|_{\cH^{\frac12}}
:= \||\Tau|^\frac12 v\|_{L^2}.
\end{equation}
Then we have 

\begin{theorem}\label{t:KS-Besov}
The following pointwise bounds hold  
\begin{equation}\label{KS}
\begin{aligned}
|\phi(x)| & \ \lesssim  \omega_0(t,x) ( \| \phi\|_{\cB^{-\frac12}_{2,\infty}} + \|L\phi \|_{\cH^{\frac12}}) , 
\\
|\Tau \phi(x)| & \ \lesssim \omega_1(t,x)
( \| \phi\|_{\cB^{-\frac12}_{2,\infty}} + \|L\phi \|_{\cH^{\frac12}}).
\end{aligned}
\end{equation}

\end{theorem}

\begin{remark}
For $x >0$ the operator $L$ is elliptic and hence better pointwise bounds are expected in this region. This justifies the result we stated above.

\end{remark}

\begin{proof}
One may think of the the bounds \eqref{KS} as fixed time bounds, without reference to any time evolution. 
In order to prove these bounds there are two possible strategies: 

\begin{enumerate}
\item Interpret $L$ as a hyperbolic operator 
in the region $x < 0$ and as an elliptic operator if $x > 0$, and directly prove suitable propagation, respectively elliptic estimates in the two regions, or 
\item consider the associated time evolution and use instead the dispersive bounds.
\end{enumerate}
Both strategies would work here. For examples where the first strategy is implemented we refer the reader to 
\cite{KIT}, \cite{IT-wp}. But here instead we will implement the second strategy, and show how our theorem can be proved using the dispersive estimates in Proposition~\ref{p:bodisp-loc}.

We think of the function $\psi$ in the theorem as the time $t$ section of a solution, still denoted by $\phi$,
for the linear ILW equation \eqref{ilw-lin}. We have the $L^2$ type conservation 
\[
\| \phi(t) \|_{\cB^{-\frac12}_{2,\infty}} = \| \phi(0) \|_{\cB^{-\frac12}_{2,\infty}},
\qquad \|L\psi (t)\|_{\cH^{\frac12}} = \|L\phi (0)\|_{\cH^{\frac12}}.
\]
To shorten the notations we normalize 
\[
\| \phi(0) \|_{\cB^{-\frac12}_{2,\infty}} + \|L\phi (0)\|_{\cH^{\frac12}} = 1.
\]

Then we use a dyadic decomposition to write
\[
\begin{aligned}
\phi(t) = & \  e^{itA} \phi(0) = e^{itA} P_{> 0} \tilde P_{>0} \phi(0) + \sum_{j \leq 0} e^{itA} P_{j} \tilde P_{j} \phi(0)
\\
= & \ K_{>0}(t) \ast \tilde P_{>0} \phi(0)+ \sum_{j \leq 0} e^{itA} K_{j} \ast  \tilde P_{j} \phi(0).
\end{aligned}
\]

We will estimate separately the terms in the above sum. We begin by estimating their initial data.
For $\tilde P_{> 0} \phi(0)$ we can commute 
\[
x \tilde P_{> 0} = P_{> 0} x + i P'_{<0} \, ,
\]
where both terms use only frequencies $\gtrsim 0$.
Then we can directly write
\[
\| \tilde P_{> 0} \phi(0)\|_{L^2} + \| \tilde P_{> 0} \phi(0)\|_{L^2}
\lesssim \| \phi(0) \|_{\cB^{-\frac12}_{2,\infty}}+ \|x \phi (0)\|_{\cH^{\frac12}} \lesssim 1 .
\]
Then we combine this with \eqref{K0} to obtain
\begin{equation}\label{k0-phi}
 |  K_{>0}(t) \ast \tilde P_{>0} \phi(0) | + |  \Tau K_{>0}(t) \ast \tilde P_{>0} \phi(0) | \lesssim   \left\{
 \begin{aligned} 
 & t^{-\frac12}  \qquad \mbox{ for } x < - \frac18 t \\
 &t^{-\frac{1}{2}} (1+|x|+t)^{-1} \mbox{ otherwise. } 
 \end{aligned}\right.
\end{equation}

Similarly for $j \leq 0$ we commute
\[
x \tilde P_{> 0} = P_{j} x + i P'_{j},
\]
where both terms use only frequencies of size $2^j$, and the symbol $\tilde P'_j$ has size $2^{-j}$.
Then we can estimate
\[
2^{-\frac{j}2}\| \tilde P_j \phi(0)\|_{L^2} + 2^{\frac{j}2}\| \tilde P_{> 0} \phi(0)\|_{L^2}
\lesssim \| \phi(0) \|_{\cB^{-\frac12}_{2,\infty}}+ \|x \phi (0)\|_{\cH^{\frac12}} \lesssim 1.
\]
Combining this with \eqref{Kj} we obtain
\begin{equation}\label{kjphi}
\begin{aligned}
 & |K_j(t,x)\ast \tilde P_j \phi(0) | + 2^{-j} |\Tau(D) K_{j}(t,x)\ast \tilde P_j \phi(0)| \lesssim 
 \\ & 
\lesssim  \left\{
 \begin{aligned} 
 & 2^{-\frac{j}2} (t+2^{-3j})^{-\frac12}  \qquad \mbox{ for }  -  2^{2j+4} t <  x < -  2^{2j-4} t \\
 & 2^{-\frac{j}2} (t+2^{-3j})^{-\frac12} (1+2^j|x|+2^{3j} t)^{-1} \mbox{ otherwise.} 
  \end{aligned}
 \right.
\end{aligned}
\end{equation}

Finally, the desired bound follows by summing up \eqref{k0-phi} and \eqref{kjphi}.  To see this, we consider first the bound for $\phi(t)$. We observe that the angular concentration regions are 
essentially disjoint, and give exactly the weights 
$\omega_1$ and $\omega_2$. Hence it remains to estimate the tail, i.e. the sum
\[
S_0 = \sum_{j\leq 0}  2^{-\frac{j}2} (t+2^{-3j})^{-\frac12} (1+2^j|x|+2^{3j} t)^{-1}.
\]
The summand is piecewise monotone in $j$ and decaying at $-\infty$, so it suffices to evaluate it at the midpoints.
We distinguish two cases:

a) If $|x| \lesssim t^{-\frac13}$ then we can discard the $|x|$, and there is only one midpoint, at $2^{3j} = t^{-1}$.
Thus we get $S_0 \approx t^{-\frac13}$.

b) If $|x| \gg t^\frac13$ then we have three midpoints,
at $2^j = |x|^{-1}$, at $2^{3j} = t^{-1}$ and finally 
at $2^{2j} = x/t$. Evaluating the summand at these three 
points we obtain the values $|x|^{-1}$, $|x|^{-1}$ respectively $t^\frac14 |x|^{-7/12}$, of which the last one is smaller and can be discarded. Then it remains to 
add up the dyadic range between the first to values of $j$, which is about $\ln (x t^{-\frac13})$.
Thus we obtain 
\[
S_0 \lesssim t^{-\frac13} \langle x t^{-\frac13} \rangle^{-1} \ln  \langle x t^{-\frac13} \rangle^{-1}.
\]

The computations are similar for the $\Tau \psi$ bound.
There we need to consider the sum
\[
S_1 = \sum_{j\leq 0}  2^{\frac{j}2} (t+2^{-3j})^{-\frac12} (1+2^j|x|+2^{3j} t)^{-1}.
\]
Case (a) is similar, but the balance changes in case 
(b), where at the midpoints we get the values
$|x|^{-2}$, $|x|^{-1} t^{-\frac13}$, respectively $t^{-\frac14} x^{-\frac54}$ of which the middle one is largest. Hence we obtain
\[
S_1 \lesssim |x|^{-1} t^{-\frac13}, 
\]
as needed.

\end{proof}

\subsection{Strichartz and bilinear \texorpdfstring{$L^2$}{} bounds}
In this section we discuss the Strichartz estimates adapted to the ILW equation, as well as corresponding 
bilinear $L^2$ bounds. These will be used in the next section in the proof of local well-posedness, so we only need them locally in time. 

The advantage in working on a bounded time interval is that there one may view the ILW flow as a perturbation of the Benjamin-Ono equation, as explained later in Section 4.2, see the equivalent form \eqref{ilw-r}.
Then the  short time Strichartz bounds for ILW are the
same as the ones for Benjamin-Ono flow
\begin{equation}\label{bo-lin-inhom}
(\partial_t + H\partial^2)\psi = f, \qquad \psi(0) = \psi_0.
\end{equation}

We define the Strichartz space $S$ associated to the $L^2$ flow by
\[
S = L^\infty_t L^2_x \cap L^4_t L^\infty_x,
\]
as well as its dual
\[
S' = L^1_t L^2_x + L^{\frac43} _t L^1_x .
\]

The Strichartz estimates in the $L^2$ setting are summarized in the following

\begin{lemma}\label{l:Str}
Assume that $\psi$ solves either \eqref{bo-lin-inhom} ILW in $[0,T] \times \mathbb{R}$. Then 
the following estimate holds.
\begin{equation}
\label{strichartz}
\| \psi\|_S \lesssim \|\psi_0 \|_{L^2} + \|f\|_{S'} .
\end{equation}
\end{lemma}

We remark that these Strichartz estimates can also be viewed as a consequence \footnote{ Except for the $L^4_tL_x^{\infty}$ bound, as the Hilbert transform is not bounded in $L^{\infty}$.}
of the similar estimates for the linear Schr\"odinger equation. This is because the two flows
agree when restricted to functions with frequency localization in $\mathbb{R}^+$.

We also remark that we have the following Besov version of the 
estimates, 
\begin{equation}
\label{strichartzB}
\| \psi\|_{\ell^2 S} \lesssim \|\psi_0 \|_{L^2} + \|f\|_{\ell^2 S'},
\end{equation}
 where 
\[
\| \psi \|_{\ell^2 S}^2 = \sum_k \| \psi _k \|_{ S}^2, \qquad  \| \psi \|_{\ell^2 S'}^2 = \sum_k \| \psi _k \|_{ S'}^2 .
\]

The last property we transfer from the linear Benjamin-Ono equation is the 
bilinear $L^2$ estimate, which is as follows:

\begin{lemma}
\label{l:bi}
Let $\psi^1$, $\psi^2$ be two solutions to the inhomogeneous Benjamin-Ono equation or the ILW equation with data  $\psi^1_0$, $\psi^2_0$ and inhomogeneous terms $f^1$ and $f^2$, in a time interval $[0,T]$. Assume 
that the sets 
\[
E_i = \{ |\xi|, \xi \in \text{supp } \hat \psi^i \}
\]
are disjoint. Then we have 
\begin{equation}
\label{bi-di}
\| \psi^1 \psi^2\|_{L^2} \lesssim \frac{1}{\text{dist}(E_1,E_2)} 
( \|\psi_0^1 \|_{L^2} + \|f^1\|_{S'}) ( \|\psi_0^2 \|_{L^2} + \|f^2\|_{S'}).
\end{equation}
\end{lemma}

These bounds also follow from the similar bounds for the Schr\"odinger equation, where only the separation 
of the supports of the Fourier transforms is required. They can be obtained in a standard manner  from the similar bound for products of solutions to the homogenous equation, for which we refer the reader to \cite{tao-m}.

One corollary of this result  applies in the case when we look at the product
of two solutions which are supported in different dyadic regions:
\begin{corollary}\label{c:bi-jk}
Assume that $\psi^1$ and $\psi^2$ as above are supported in dyadic regions $\vert \xi\vert \approx 2^j$ and $\vert \xi\vert \approx 2^k$, $\vert j-k\vert >2$, then
\begin{equation}
\label{bi}
\| \psi^1 \psi^2\|_{L^2} \lesssim 2^{-\frac{\max\left\lbrace j,k \right\rbrace  }{2}}
( \|\psi_0^1 \|_{L^2} + \|f^1\|_{S'}) ( \|\psi_0^2 \|_{L^2} + \|f^2\|_{S'}).
\end{equation}
\end{corollary}

Another useful  case is when we look at the product
of two solutions which are supported in the same  dyadic region, but with frequency separation:
\begin{corollary}\label{c:bi-kk}
Assume that $\psi^1$ and $\psi^2$ as above are supported in the dyadic region 
$\vert \xi\vert \approx 2^k$, but have $O(2^k)$ frequency separation between their supports.
Then 
\begin{equation}
\label{bi-kk}
\| \psi^1 \psi^2\|_{L^2} \lesssim 2^{-\frac{k}{2}}
( \|\psi_0^1 \|_{L^2} + \|f^1\|_{S'}) ( \|\psi_0^2 \|_{L^2} + \|f^2\|_{S'}).
\end{equation}
\end{corollary}

\section{Normal form analysis and local well-posedness}


In this section we prove the local well-posedness result in Theorem~\ref{thm:lwp}. 
Since the $L^2$ norm of the solution is conserved, this in turn implies global well-posedness. For local well-posedness, it is the high frequency behavior of the solutions which matters most, and this 
is where the ILW equation asymptotically coincides with the Benjamin-Ono equation. Because of this, we will interpret the ILW 
flow as a perturbation of the Benjamin-Ono flow.  But both flows are quasilinear, so a direct perturbative approach is out of the question.

Nevertheless, the ideas we will use in the proof are directly inspired by the work of Ifrim-Tataru \cite{IT-BO} where, among other results, a low regularity global well-posedness result was proved for the Benjamin-Ono equation; this was done without relying on any complete integrability specific tools.

The first step in the proof is to reduce the problem to the small $L^2$ data case. In the Benjamin-Ono case, this reduction is achieved directly by scaling. Here, rescaling large data
int small data also has the effect of changing the equation, precisely the parameter $\delta$, from $\delta = 1$ into a large
parameter. Once this is achieved, in order to improve the analogy with Benjamin-Ono we also eliminate the transport term via a linear change of coordinates, see Remark~\ref{r:transport}. 

After the above transformations, we have reduced the problem to proving local well-posedness for the problem \eqref{ilw-} with $\delta \geq 1$ where we assume that the initial data satisfies
\begin{equation}
\label{small id}
\Vert \phi(0)\Vert_{L^2_x}\leq \epsilon,
\end{equation}
where $\epsilon$ is a small universal constant that does not depend on the choice of $\delta$. One advantage of working with small initial data is that now we can fix the time interval where 
we prove the local well-posedness result to $I = [-1,1]$.

Following the succession of steps in \cite{IT-BO}, we split the proof as follows:

 \begin{enumerate}[label=\roman*)]
     \item We establish apriori $L^2$ bounds for regular $H^{\frac{3}{2}}$ solutions for the problem \eqref{ilw-}.
 \item
 We prove $\dot H^{-\frac12}$ bounds for the linearized equation.
 \item 
 We combine the first two steps to construct $L^2$ solutions as limits of $H^\frac32$ solutions and conclude the proof of the theorem. 
 \end{enumerate}
Once the first two steps are carried out, the last step follows from
general principles, see the similar argument in \cite{IT-BO} and also the quasilinear well-posedness primer \cite{IT-primer}. So we will focus on the first two steps, for which we state the results in the following two theorems. We begin with the a-priori bounds for the solutions.

\begin{theorem} \label{t:apriori} 
Let $\phi$ be an $H^{\frac{3}{2}}_x$ solution to \eqref{ilw} with small initial data as in \eqref{small id}. Let $\left\lbrace
    c_k\right\rbrace _{k=0}^{\infty}\in l^2 $ so that $\epsilon c_k$
  is a frequency envelope for the initial $\phi(0)$ in $L^2$.
 Then we have the Strichartz bounds
\begin{equation}\label{u-small}
 \Vert \phi_k \Vert_{S^{0}([-1,1] \times \mathbb{R})} \lesssim \epsilon c_k,
\end{equation}
as well as the bilinear bounds
 \begin{equation}
 \label{bilinear-small}
  \Vert \phi_j \cdot \phi_k\Vert_{L^2}\lesssim 2^{-\frac{\max \left\lbrace j,k\right\rbrace}{2} }\epsilon^2 c_k\, c_j,  \qquad j\neq k.
   \end{equation}
   \end{theorem} 

These bounds are phrased and proved here for $H^{\frac{3}{2}}$ solutions, but at the conclusion of the argument they also follow for $L^2$ solutions.

Next we state the bounds for the linearized equation.

\begin{theorem} \label{t:linearized} 
Let $\phi$ be an $H^{\frac{3}{2}}_x$ solution to \eqref{ilw} with small initial data as in \eqref{small id}. Then the linearized equation around $\phi$ is well-posed in $H^{-\frac12}$, with a uniform bound
\begin{equation}
\|v\|_{S^{-\frac12}([0,1] \times \mathbb R)} \lesssim \|v(0)\|_{H^{-\frac12}}.    
\end{equation}
\end{theorem}

Here it is essential to study the linearization at lower regularity,
and the space $H^{-\frac12}$ is in some sense best suited for this purpose. Indeed, studying the linearized problem in $L^2$  would yield Lipschitz dependence in $L^2$ for the solution to data map, which is known to be false for the Benjamin-Ono equation and also for its perturbations as it is the case here. 
   
 The implicit constants in both theorems are universal, and in particular  do not depend on the $H^{\frac3{2}}_x$ norm of the initial data $\phi(0)$.  

 We remark that as part of the proof of Theorem~\ref{t:linearized}
 we also prove bilinear $L^2$ bounds, which are stated later. However these bounds are not needed in order to complete the proof of the main well-posedness result in Theorem~\ref{thm:lwp}.

 The rest of the section is devoted to the proof of these two theorems. As a general principle, we note that a standard iteration method will not work,   because the linear part of the Intermediate Long Wave equation does not have enough smoothing to compensate for the derivative in the nonlinearity. In fact, the ILW equation should be seen as a perturbation of the Benjamin-Ono  equation, as we will emphasize later in the proof. To resolve the difficulty arising from the lack of a standard iteration method we use a succesion of  ideas as in the work of Ifrim-Tataru \cite{IT-BO}, and which are related to   the normal form method, first introduced by Shatah in \cite{shatah} in the context of dispersive PDEs. The main principle in the
   normal form method is to apply a quadratic correction to the
   unknown in order to replace a nonresonant quadratic nonlinearity by
   a milder cubic nonlinearity. Unfortunately this method does not
   apply directly here, because some terms in the quadratic correction
   are unbounded, and so are some of the cubic terms generated by the
   correction. To bypass this issue here we develop a more favorable
   implementation of normal form analysis. This is carried out in two
   steps:
   \begin{itemize}
   \item  a partial normal form transformation which is bounded and removes some of the quadratic nonlinearity
  \item a conjugation via a suitable exponential (also called gauge transform, \cite{tao}) which removes in a bounded way the remaining part of the quadratic nonlinearity.
   \end{itemize}
 This will transform the Intermediate Long Wave equation \eqref{ilw} into an equation where the quadratic terms have been removed  and replaced by cubic perturbative terms. A similar approach applies in the study of the linearized equation.

\subsection{The quadratic normal form analysis} In this subsection we rewrite the equation \eqref{ilw-} as a perturbation of the Benjamin-Ono equation, and carry out the analysis as in \cite{IT-BO}. The equivalent form of the  equation is
\begin{equation}
\label{ilw-r}
\phi_t +H\partial^2_x \phi =\phi\phi_x +\mathcal{P}\partial_x^2\phi ,\quad  \mathcal{P}:=\left( H-\Tau^{-1}_\delta\right) ,
\end{equation}
where $H$ is the Hilbert transform and  $\mathcal{P}:=H-\Tau_\delta^{-1}$ is a multiplier with a bounded  symbol that decays exponentially at high frequency.  This will allow us to treat it in a perturbative manner in order to reach our conclusion. 

Before going further, we emphasize that by a \emph{normal form} we
refer to any type of transformation which removes nonresonant
quadratic terms; all such transformations are uniquely determined up
to quadratic terms.

We recall that the quadratic  normal form transformation for the Benjamin-Ono equation  was formally derived in  \cite{IT-BO}.
Even though  this was not used directly, portions 
of it were used in order  to remove certain ranges of frequency interactions from the quadratic nonlinearity.

In the study of equation \eqref{ilw-r}, we will use the same normal form as in the Benjamin-Ono case, which for convenience, we recall here:

\begin{proposition} The formal quadratic normal form transformation associated to the Benja\-min-Ono equation \eqref{BO}  is given by 
\begin{equation}
\label{bo nft}
\tilde{\phi}=\phi -\frac{1}{4}H\phi \cdot \partial^{-1}_x\phi  -\frac{1}{4}H\left( \phi \cdot \partial^{-1}_x\phi \right) .
\end{equation}
\end{proposition}
Note that at low frequencies \eqref{bo nft} is not invertible, which
tends to be a problem if one wants to apply the normal form
transformation directly.
For the proof we direct the reader to \cite{IT-BO}.

 \subsection{A modified normal form analysis}
\label{s:local}

We begin by writing the Intermediate Long Wave  equation \eqref{ilw-r} in a
paradifferential form, i.e., we localize ourselves at a frequency $2^k$,
and then project the equation either onto negative or positive
frequencies:
\[
(\partial_t \mp  i\partial^{2}_{x} )\phi_k^{\pm}=P_k^{\pm} (\phi \cdot \phi_{x}) + P_k^{\pm} \left[ \mathcal{P}\partial_x^2\phi \right].
\]
Since $\phi$ is real, $\phi^- $ is the complex conjugate of $\phi^+$ so it suffices to work with the latter. However, the case $k = 0$ is special  in the same way  it was in the Benjamin-Ono  case.  We will address it separately at each step of our analysis. 

Thus, the ILW  equation for the positive frequency Littlewood-Paley components $\phi^{+}_k$ is
\begin{equation}
\label{eq-loc}
\begin{aligned}
&\left( i\partial_{t}  + \partial^{2}_{x} \right)  \phi^+_{k}=iP_k ^{+}(\phi \cdot \phi_{x}) +iP_k^{+} \left[  \mathcal{P}\partial_x^2\phi \right] .
\end{aligned}
\end{equation}

Heuristically, the worst term in $P_k^+(\phi \cdot \phi_{x})$ occurs
when $ \phi _x$ is at high frequency and $\phi$ is at low
frequency. We can approximate $P_k^+ (\phi \cdot \phi_x)$, by its
leading paradifferential component $\phi_{<k} \cdot \partial_x
\phi_{k}^+$; the remaining part of the nonlinearity will be harmless.
More explicitly we can eliminate it by means of a bounded normal form
transformation.

We will extract out the main term
$i\phi_{<k}\cdot \partial_{x}\phi_k^+$ from the right hand side
nonlinearity and move it to the left, obtaining
\begin{equation}
\label{eq-lin-k}
\left( i\partial_t+\partial_x^2 -i\phi_{<k}\cdot \partial_{x}\right)  \phi^+_k =i P_{k}^{+}\left( \phi_{\geq k} \cdot \phi _{x}\right) +i\left[ P_{k}^{+}\, , \, \phi_{<k} \right]\phi_x
+ iP_k^{+} \left[  \mathcal{P}\partial^2_x\phi \right].
\end{equation}
For reasons which will become apparent later on, when we do the
exponential conjugation, it is convenient to add an additional lower
order term on the left hand side (and thus also on the right).
Denoting by $A^{k,+}_{IT-BO} $ the operator
\begin{equation}
\label{op Abo}
A^{k, +}_{IT-BO}:= i\partial_t +\partial^2 _x-i\phi_{<k}\cdot \partial_{x}  +\frac{1}{2}  \left( H+i \right)\partial_x\phi _{<k}
\end{equation} 
we rewrite the equation \eqref{eq-lin-k} in the form
\begin{equation}
\label{eq-lin-k+}
A^{k,+}_{IT-BO} \ \phi^+_k =i P_{k}^{+}\left( \phi_{\geq k} \cdot \phi _{x}\right) +i\left[ P_{k}^{+}\, , \, \phi_{<k} \right]\phi_x +\frac{1}{2}  \left( H+i \right)\partial_x\phi _{<k} \cdot \phi^+_k +iP_k^{+} \left[  \mathcal{P}\partial^2_x\phi \right].
\end{equation}
Note the key property that the operator $A^{k,+}_{IT-BO} $ is symmetric, which in particular tells us that the $L^2$ norm of the solution for the associated linear equation is conserved in the
corresponding linear evolution.

 The case $k=0$  needs a separate discussion. There one does not  need to work with the problem in a paradifferential form, in order to avoid the operator $P_0^+$ which does not have a smooth symbol.
Thus we will work with the equation
\begin{equation}
\label{eq-lin-0}
(\partial_t + H \partial_x^2)\phi_0 = P_0( \phi_{0} \phi_x) + P_0( \phi_{>0} \phi_x) +P_0 (\mathcal{P}\partial_x^2),
\end{equation}
where the first term on the right is purely a low frequency term and will play only a perturbative role.

The next step is to eliminate the nonresonant quadratic terms on the right hand side of
\eqref{eq-lin-k+} using a normal form transformation
 \begin{equation}
\label{partial nft}
\begin{aligned}
 \tilde{\phi}_k^+:= \phi^+_k +B_{k}(\phi, \phi).
\end{aligned}
\end{equation} 
Such a transformation is easily computed and formally is given by the expression
\begin{equation}
\label{bilinear nft}
\begin{aligned}
B_{k}(\phi, \phi) =&\frac{1}{2}HP_{k}^{+}
\phi \cdot \partial_{x}^{-1}P_{<k}\phi-\frac{1}{4}P_{k} ^{+}\left( H\phi\cdot \partial_{x}^{-1}\phi\right)
 -\frac{1}{4}P_{k}^{+}H\left( \phi \cdot \partial^{-1}_x\phi\right).
\end{aligned}   
\end{equation}
One can view this as a subset of the normal form transformation
computed for the full equation, see \eqref{bo nft}. Unfortunately, as
written, the terms in this expression are not well defined because
$\partial^{-1}_x \phi$ is only defined modulo constants. To avoid this
problem we separate the low-high interactions, which yields a well
defined commutator, and we rewrite $B_{k}(\phi, \phi)$  as
\begin{equation}
\label{commutator B_k}
B_{k}(\phi, \phi)=-\frac{1}{2}\left[ P^+_kH\, , \, \partial^{-1}_x \phi_{<k} \right]\phi -\frac{1}{4}P^+_k \left( H\phi \cdot \partial_x^{-1}\phi_{\geq k}\right)  
-\frac{1}{4}P^+_k H \left(\phi \cdot \partial_x^{-1}\phi_{\geq k}\right).
\end{equation}

Replacing $\phi^+_k$ with $\tilde{\phi}^+_k $ removes all the quadratic
terms on the right and leaves us with an equation of the form
 \begin{equation}
 \label{eq after 1nft}
  A^{k,+}_{IT-BO} \, \tilde{\phi}^{+}_k = R_{k}^1(\phi)+R_{k}^{2}(\phi, \phi) +
  Q^3_k(\phi, \phi, \phi ),
 \end{equation}
 where $Q^{3}_k(\phi, \phi, \phi )$ contains only cubic  terms in $\phi$, while the term $R_{k}^{2}(\phi, \phi)$ contains only quadratic terms. 
 The expression of  $Q^{3}_k(\phi, \phi, \phi)$ is the same as in the Benjamin-Ono case, and is 
 given  in greater detail and examined  in \cite{IT-BO}, Lemma~$4.3$. On the other hand, the last term on the right in \eqref{ilw-r} produces some new contributions, a quadratic and a linear one, represented by the last two terms in the above formula.
 For clarity, here and below we use the letter $R$ for these new
 contributions. The linear term $R_{k}^1(\phi)$ is the same as 
 in \eqref{eq-lin-k+}, namely 
 \[
R_{k}^1(\phi) = iP_k^{+} \left[  \mathcal{P}\partial^2_x\phi \right].
 \]
The expression and the bounds on $R^{2}_k (\phi, \phi)$  will be discussed later in Lemma~\ref{l:Rk}.

 We return to the case $k=0$, where we note that  here the first normal form transformation does not eliminate the low-low frequency interactions, and our intermediate  equation has the form
 \begin{equation}
 \label{eq after 1nft-0}
(i\partial_t +\partial^2_x) \, \tilde{\phi}_0 =R^1_0 (\phi)+ R^2_0 (\phi, \phi) + Q^2_0 (\phi, \phi) +Q^{3}_0(\phi, \phi, \phi ),
 \end{equation}
where $Q^2_0$ contains all the low-low frequency interactions
\[
Q^2_0 (\phi, \phi):=P_0 \left(  \phi_{0}\cdot \phi_x\right). 
\]

The second stage in our normal form analysis is to perform a second bounded normal form transformation that will remove the paradifferential terms in
the left hand side of \eqref{eq after 1nft}; this will be a renormalization, following the idea introduced by Tao \cite{tao}. To achieve this we  introduce and initialize the spatial primitive $\Phi(t, x)$ of $\phi (t,x)$, exactly as in \cite{tao}.
 It turns out that $\Phi (t,x)$ is necessarily a real valued function that  solves the equation
\begin{equation}
\label{Phi}
\Phi_{t}+H\Phi_{xx}=\Phi_{x}^2+ \mathcal{P}\Phi_{xx},
\end{equation}
which holds globally in time and space. Here, the constants are fixed by imposing the initial condition  $\Phi(0,0)=0$. Thus,  
\begin{equation}
\label{antiderivative}
\Phi_x (t,x)= \frac{1}{2}\phi(t,x).
\end{equation}
The idea in \cite{tao} was that in order to get bounds on $\phi$ it suffices to obtain appropriate bounds 
 on $\Phi(t,x)$ which are one higher degree of  regularity, as \eqref{antiderivative} suggests. 
Here we instead use $\Phi$ merely in an auxiliary role, in order to define the second normal form transformation. This is
 \begin{equation}
 \label{conjugation}
 \displaystyle{\psi_k^+:=\tilde{\phi}_k^{+} \cdot    e^{-i\Phi_{<k}}}.
 \end{equation}

 The transformation \eqref{conjugation} is akin to a Cole-Hopf
 transformation, and expanding it up to quadratic terms, one observes
 that the expression obtained works as a normal form transformation,
 i.e., it removes the paradifferential quadratic terms.  The
 difference is that the exponential will be a bounded transformation,
 whereas the corresponding quadratic normal form is not.  One also
 sees the difference reflected at the level of cubic or higher order
 terms obtained after implementing these transformation (obviously
 they will differ).

 The case $k = 0$ is special here as well, in that this  renormalization  step is  not needed. There we
simply set $\psi_0 =  \tilde \phi_0$, and use the equation \eqref{eq after 1nft-0}.

By applying this \emph{Cole-Hopf type} transformation, we  rewrite  the equation \eqref{eq after 1nft}
as a nonlinear Schr\"odinger equation for our final normal form variable $\psi_k^+$
 \begin{equation}
 \label{conjugare}
 \begin{aligned}
(i\partial_t +\partial^2_x)\, \psi_k^+ = 
[\tilde{R}^{1}_k(\phi )+\tilde{R}^{2}_k(\phi ,\phi )
+\tilde{Q}_k^{3}(\phi, \phi, \phi) + \tilde{Q}_k^{4}(\phi, \phi, \phi,\phi)] e^{-i\Phi_{<k}}.
 \end{aligned}
 \end{equation}
Using the same style for notations as before, the terms 
$\tilde{Q}_k^{3}$ and $\tilde{Q}_k^{4}$ are the same as in \cite{IT-BO}, whereas the terms $\tilde R^1_k$, $\tilde R^2_k$  are the new ones. Here $\tilde R^1_k = R^1_k$, while $\tilde R^2_k$ is also computed in
 Lemma~\ref{l:Rk}.
 
 \begin{lemma}\label{l:Rk}  The quadratic, cubic and quartic expressions $R_{k}^2(\phi, \phi)$, $\tilde{R}^2_k (\phi, \phi)$, $\tilde{Q}_k^3$ and $\tilde{Q}^4_k$ are translation invariant multilinear forms of the type
 \begin{equation}
 \begin{aligned}
\label{Rk}
&R^2_k(\phi, \phi)=L_k (\phi, \mathcal P \phi_x)+L_k(H\phi,  \mathcal P \phi_x)  +L_k(\phi,  H \mathcal P \phi_x),     \\
&\tilde{R}^2_k(\phi, \phi)=
L_k (\phi, \mathcal P \phi_x)+L_k(H\phi,  \mathcal P \phi_x)  +L_k(\phi,  H \mathcal P \phi_x),\\
&\tilde{Q}^{3}_k(\phi, \phi, \phi)=L_k (\phi, \phi, \phi) +L_k(H\phi, \phi, \phi),\\
& \tilde{Q}^4_k (\phi, \phi, \phi, \phi)= L_k (\phi, \phi, \phi ,\phi) +L_k(H\phi, \phi, \phi, \phi).
 \end{aligned}
  \end{equation}
 \end{lemma}
 \begin{proof}
 We begin by computing the expression for $R_k^{2}$,
 \[
\begin{aligned}
R_{k}^{2}(\phi, \phi):=&+\frac{1}{2}iHP_{k}^{+} \left\{  \mathcal{P}\phi_{xx}\right\} \cdot \partial_{x}^{-1}P_{<k}\phi +\frac{1}{2}iHP_{k}^{+} \phi \cdot P_{<k}\left\{  \mathcal{P}\phi_{x}\right\}\\
&-\frac{1}{4}iP_{k} ^{+}\left( H\left\{  \mathcal{P}\phi_{xx}\right\}\cdot \partial_{x}^{-1}\phi\right)-\frac{1}{4}iP_{k} ^{+}\left( H\phi\cdot \left\{  \mathcal{P} \phi_{x}\right\}\right)\\
&-\frac{1}{4}iP_{k}^{+}H\left( \left\{  \mathcal{P}\phi_{xx}\right\} \cdot \partial^{-1}_x\phi\right) -\frac{1}{4}iP_{k}^{+}H\left( \phi \cdot \left\{  \mathcal{P}\phi_{x}\right\}\right).
\end{aligned}
\]
 To avoid inverse derivatives at low frequency we rewrite it in a commutator fashion using
\[
\begin{aligned}
&\left[ \partial_x^{-1}\phi_{<k}\, ,\, HP^+_k\right]  \mathcal{P}\phi_{xx}=HP^+_k \mathcal{P}\phi_{xx}\cdot \partial_x^{-1}\phi_{<k}-
HP^+_k\left[  \mathcal{P}\phi_{xx}\cdot \partial^{-1}_{x}\phi_{<k}\right],\\ 
&\left[ \partial^{-1}_x\phi_{<k}\, ,\, P^+_k \right]H \mathcal{P} \phi_{xx} =P^+_kH \mathcal{P}\phi_{xx}\cdot \partial^{-1}_x\phi_{<k}
-P^+_k\left[ H \mathcal{P}\phi_{xx}\cdot \partial^{-1}_x\phi_{<k}\right],
\end{aligned}
\]
so that
\[
\begin{aligned}
R_{k}^2(\phi, \phi)
&=\frac{1}{4}i\left[ \partial_x^{-1}\phi_{<k}\, ,\, HP^+_k\right]  \mathcal{P}\phi_{xx} +\frac{1}{4}i\left[ \partial^{-1}_x\phi_{<k}\, ,\, P^+_k \right]H \mathcal{P}\phi_{xx}\\
&-\frac{1}{4}iP_{k} ^{+}\left( H \mathcal{P}\phi_{xx}\cdot \partial_{x}^{-1}\phi_{\geq k}\right)-\frac{1}{4}iP_{k} ^{+}\left( H\phi\cdot   \mathcal{P}\phi_{x}\right)\\
&-\frac{1}{4}iP_{k}^{+}H\left(   \mathcal{P}\phi_{xx} \cdot \partial^{-1}_x\phi_{\geq k}\right) -\frac{1}{4}iP_{k}^{+}H\left( \phi \cdot  \mathcal{P}\phi_{x}\right)\\
&+\frac{1}{2}iHP_{k}^{+} \phi \cdot P_{<k} \mathcal{P}\phi_{x} .\\
\end{aligned}
\]
For the first two terms we use Lemma~\ref{commutator} and eliminate the low frequency antiderivative term. We further apply Lemma~\ref{l:hilbert-tilbert} to cancel the remaining derivative in the other expressions. For instance,  consider the first commutator
\[  
\begin{aligned}
\left[ \partial^{-1}\phi_{<k}, HP^{+}_k\right]  \mathcal{P} \phi_{xx}&=\left[ \partial^{-1}\phi_{<k}, HP^{+}_k\right] \tilde{P}_k \mathcal{P}\phi_{xx}\\
&=L_k\left(\phi_{<k}, 2^{-k}\tilde{P}_k  \mathcal{P}\phi_{xx}\right)\\
&=L_k\left( \phi_{<k}, \mathcal P \phi_x \right).
\end{aligned}
\]
  For the remaining terms we split the unlocalized $\phi$ factor into $\phi_{\geq k} +\phi_{<k}$. The contribution of $\phi_{<k}$ is as before, while the remaining bilinear term in $\phi_{\geq k}$ the frequencies of the two inputs must be balanced at some frequency $2^j$ where $j$ ranges in the region $j\geq k$.  The terms  in the expression of $R_k^2$ are all similar.

We can also consider, for example, the third term of $R_k^2$  
\[
P^+_k\left( H  \mathcal{P}\phi_{xx} \cdot \partial_x ^{-1}\phi_{\geq k}\right) = P^+_k\partial_x \left( H \mathcal{P} \phi_x \cdot \partial_x^{-1}\phi_{\geq k}\right)-P^+_k\left( H \mathcal{P} \phi_x \cdot \phi_{\geq k}\right). 
\]
The first derivatives in the first term above yields a $2^k$ factor and we can use the following observation
\begin{equation}
\label{high-derivative}
\partial^{-1}_x\phi_{\geq k}=2^{-k}L(\phi)
\end{equation}
to cancel the $2^k$ factor. It follows that indeed
\[
R^2_k(\phi, \phi)=L_k(H\phi,  \mathcal P\phi_x)+L_{k}(\phi, \mathcal P\phi_x)  +  L_k(\phi,  H \mathcal P\phi_x)
\]
as needed.

The expression of $\tilde{R}^2_k$ is given by
\[
\tilde{R}^2_k(\phi, \phi):= \frac12 \phi^+_k\cdot  \mathcal{P}\phi_{x,<k}
+  R^2_k(\phi, \phi),
\]
and the first  term is already in the desired form.

For $\tilde{Q}_k^{3}$ and $\tilde{Q}_k^4$ we redirect the reader to the work in \cite{IT-BO}.
For completeness, their expressions are given below
\[
\tilde{Q}_k^{3}(\phi, \phi, \phi):=\left( Q^3_{k}(\phi, \phi, \phi) +\frac14 \phi^+_k\cdot P_{<k}\left(\phi^2\right) +\frac12 B_k( \phi, \phi)\cdot  \mathcal{P}\phi_{x, <k}-\frac14\phi^+_k\cdot \left( \phi_{<k}\right)^2\right) e^{-i\Phi_{<k}},
\]
respectively
\[
\tilde{Q}^4_k (\phi, \phi, \phi, \phi)=\frac14 B_{k}(\phi, \phi) \cdot \left\lbrace  2P_{<k}(\phi ^2) -\left( P_{<k}\phi\right)^2   \right\rbrace .
\]
 \end{proof}

 \subsection{The bootstrap argument} We return to the proof of Theorem \ref{t:apriori}. The method we will use is a standard continuity argument based on the $H^{\frac{3}{2}}$ global well-posedness theory. For $0<t_0\leq 1$ we introduce the following norm we want to track
 \[
 M(t_0):= \sup_{k} c_k^{-2} \,  \|  P_k \phi\|_{S^0[0,t_0]}^2 + \sup_{j\neq k\in \mathbb{N}} \sup_{y\in \mathbb{R}} c^{-1}_j \cdot c_k^{-1} \cdot \| \phi_j \cdot T_y \phi_k \|_{L^2\left[ 0, t_0\right] }. 
\]
Here, in the first term, we have the Strichartz dyadic norm associated to the $L^2$ flow together with its corresponding  frequency envelope. In the second term, the role of the condition $j\neq k$ is to
insure that $\phi_j$ and $\phi_k$ have $O(2^{\max \left\lbrace j,k\right\rbrace })$ separated frequency localizations.  However,
by a slight abuse of notation, we also allow bilinear expressions of
the form $ P_{k}^1\phi \cdot P_k^2\phi$, where $P_{k}^1$ and $P_{k}^2$
are both projectors at frequency $2^k$ but with at least $2^{k-4}$
separation between the \textbf{absolute values} of the frequencies in
their support. We also want to explain the role of the  translation operator $T_y$. This is needed in order for us to be able to use thee bilinear bounds in estimating multilinear  $L$ type expressions.\\

We seek to show that 
\[
M(t)\lesssim \epsilon^2.
\]
 As $\phi$ is an $H^{\frac{3}{2}}$ solution, it is easy to see that $M(t)$
 is continuous as a function of $t$, and 
 \[
 \lim_{t\searrow 0}M(t) \lesssim \epsilon^2 .
 \]

 Thus, by a continuity argument 
it suffices to make the bootstrap assumption 
\[
M(t_0) \leq C^2 \epsilon^2 
\]
and then show that  
\begin{equation}\label{want}
 M(t_0) \lesssim  \epsilon^2 + C^4 \epsilon^4.   
\end{equation}

This bound gives the desired result provided that $C$ is is large enough (independent of
$\epsilon$) and $\epsilon$ is sufficiently small (depending on $C$).
Without any restriction in generality we can assume that $C\epsilon  \leq 1$.
From here on $t_0 \in (0,t]$ is fixed and not needed in the argument, so we drop it from the notations.
 
Given our bootstrap assumption, we have the starting estimates
\begin{equation}
\label{boot1}
  \Vert \phi_k\Vert_{S^0}\lesssim C \epsilon c_k ,
   \end{equation}
and
\begin{equation}
  \label{boot2}
  \Vert \phi_j \cdot T_y \phi_k\Vert_{L^2}\lesssim 2^{-\frac{\max \left\lbrace j,k\right\rbrace}{2} }C^2 \epsilon^2 c_j c_k,  \qquad j\neq k, \qquad y \in \mathbb{R},
\end{equation}
where in the bilinear case, as discussed above, we also allow $j=k$ provided the two localization multipliers
are at least $2^{k-4}$ separated. This separation threshold is fixed once and for all. On the other hand, when 
we prove that the bilinear estimates hold, no such sharp threshold is needed.

Our strategy will be to establish these bounds for the normal form variables $\psi_k$, and then to transfer
them to the original solution $\phi$ by inverting the normal form transformations and estimating errors.

We need to obtain bounds for the normal form variables $\psi_k^+$. 
We first consider the initial data, for which we have
\begin{lemma}
\label{l:invertibila}
Assume that the initial data $\phi_0$ satisfies \eqref{small id},
with frequency envelope $\epsilon c_k$.  Then we have 
\begin{equation}
\| \psi_k^+(0)\|_{L^2} \lesssim c_k\epsilon.
\end{equation}
\end{lemma}

This result was already proved in \cite{IT-BO}, therefore we do not include its proof. 

Next we need obtain bounds for the  right hand side in the  Schrodinger equation \eqref{conjugare} for $\psi_k^+$ in  $L^1 L^2$, which will differ from the bounds in \cite{IT-BO} by some additional  terms which will turn out to be perturbative; the bounds are obtained in the following lemma

\begin{lemma}\label{l:source}
Assume that the initial data $\phi_0$ satisfies \eqref{small id} and that the bootstrap bounds \eqref{boot1} and \eqref{boot2} hold
in $[0,1]$. Then we have 
\begin{equation}
\label{perturbative-new}
\| \tilde R^1_k \|_{L^1 L^2} \lesssim \epsilon c_k, \qquad
\| \tilde R_k^2 \|_{L^1 L^2} \lesssim C^2 \epsilon^2 c_k,
\end{equation}
respectively
\begin{equation}
\label{perturbative}
\| \tilde Q_k^3 \|_{L^1 L^2} \lesssim C^3 \epsilon^3c_k,
\qquad 
\| \tilde Q_k^4 \|_{L^1 L^2} \lesssim C^4 \epsilon^4 c_k.
\end{equation}
Further, in the case $k=0$ we also have
\begin{equation}\label{perturbative+}
  \| Q_0^2 \|_{L^1 L^2} \lesssim C^2 \epsilon^2. 
\end{equation}

\end{lemma}
\begin{proof}
The estimates \eqref{perturbative} and \eqref{perturbative+}
were proved in \cite{IT-BO}. It remains to consider 
\eqref{perturbative-new}. In the first estimate it is important that we do not have $C$. This is easily achieved using the 
mass conservation for ILW, which yields 
\[
\| \phi \|_{L^\infty L^2} \lesssim \epsilon.
\]
Then we simply use the symbol decay for $\mathcal P$ at high frequencies, together with H\"older's inequality in time. 

The bound for $\tilde R_k^2$ on the other hand uses Lemma~\ref{l:Rk}. Considering for instance the term $L_k(\phi,\mathcal P\phi_x)$,  we use the Littlewood-Paley trichotomy 
to write it as 
\[
L_k(\phi,\mathcal P\phi_x) = L(\phi_{<k} ,\mathcal P\phi_{k,x})
+ L(\phi_k,\mathcal P\phi_{<k,x}) +
\sum_{j \geq k} L(\phi_j,\mathcal P\phi_{j,x}).
\]
Then for the first two terms we use the bilinear $L^2$ bootstrap bound \eqref{boot2}, while for the last term we use twice the Strichartz bound \eqref{boot1}. In both cases we gain from the rapid decay of the symbol of $\mathcal P$ at high frequency, see Lemma~\ref{l:hilbert-tilbert}.

\end{proof}

Using the bounds in Lemmas~\ref{l:invertibila},\ref{l:source}, together with the bounds in Lemmas~\ref{l:Str},\ref{l:bi},
we obtain linear and bilinear bounds for $\psi_k^+$, namely
\begin{equation}
\label{psi1}
  \Vert \psi_k^+\Vert_{S^0}\lesssim (\epsilon + C^2 \epsilon^2) c_k ,
   \end{equation}
and
\begin{equation}
  \label{psi2}
  \Vert \psi_j^+ \cdot T_y \psi_k^+\Vert_{L^2}\lesssim 2^{-\frac{\max \left\lbrace j,k\right\rbrace}{2} }(\epsilon^2 + C^4 \epsilon^4) c_j c_k,  \qquad j\neq k, \qquad y \in \mathbb{R}.
\end{equation}
Here we can use the assumption $C \epsilon \leq 1$ to drop the 
$C$ factors. 

It remains to transfer these bounds to $\phi_k$. But this is done exactly as in \cite{IT-BO}. This concludes the proof of \eqref{want}, and thus the proof of Theorem~\ref{t:apriori}.

\subsection{The linearized equation}\label{s:linearized}
In this subsection we consider the linearized ILW equation equation,
\begin{equation}
\label{ilw-linearized}
\left\{
\begin{aligned}
& (\partial_t  +\Tau_{\delta}^{-1}\partial_x^2)v=\partial_x (\phi v) \\
&v(0)=v_0,
\end{aligned}
\right.
\end{equation}
where $v$ is the linearized variable.
This is equivalent to 
\begin{equation}
\label{ilw-linearized-P}
\left\{
\begin{aligned}
&(\partial_t +H\partial_x^2)v=\partial_x (\phi v)+\mathcal{P}\partial_x^2v\\
&v(0)=v_0.
\end{aligned}
\right.
\end{equation}

The form \eqref{ilw-linearized-P} of the linearized equation
suggests that the analysis should be performed in the same way as in the Benjamin-Ono equation, where the last term on the right can be treated  perturbatively, as the operator $\mathcal P \partial_x^2$ is bounded in all Sobolev spaces, and in particular 
in $H^{-\frac12}$. 

The proof of this well-posedness result follows the same path as the proof of the apriori bounds in the previous subsection.
A full argument is given in \cite{IT-BO} for the Benjamin-Ono case.
Some modifications are needed here, but these are similar to the modifications already given in the previous subsection for the full equation. For this reason we omit here the full details and instead
we limit ourselves to a brief outline of the argument. The steps are as follows:
\bigskip

\textbf{ Step 1.} The well-posedness result is phrased in a stronger, frequency envelope based form, where both Strichartz estimates and unbalance bilinear $L^2$ bounds are included.

\bigskip

\textbf{ Step 2.}  To facilitate the proof of the Strichartz and 
bilinear $L^2$ bounds, the proof is phrased as a bootstrap argument.

\bigskip

\textbf{ Step 3.} The linearized equation is turned into a family
of frequency localized equations for the functions $v_k^\pm = P_k^\pm v$, akin to the equations \eqref{eq-loc}.

\bigskip

\textbf{ Step 4.} The dyadic equations for $v_k^{\pm}$ are renormalized using a two step normal form transformation, which is bounded in $\dot H^{-\frac12}$. 
 \begin{itemize}
   \item  a bounded quadratic partial normal form transformation;
  \item an exponential conjugation via a suitable exponential. 
   \end{itemize}

\bigskip

\textbf{ Step 5.} The Strichartz and bilinear $L^2$ bounds 
are proved for the renormalized variables, call them $w_k^{\pm}$,
using the bootstrap assumptions, and then transferred back to $v_k^\pm$.

\bigskip

Compared to the proof of the apriori bounds in the  previous subsection, there is one difference, namely that one source term in the dyadic paradifferential equations, $P_k^\pm \partial_x(v_0 \phi)$, cannot be renormalized but is instead estimated perturbatively.

Compared to the proof of the similar result for Benjamin-Ono equation, the difference also occurs in the renormalization step, namely when using the $\phi$ equation. This yields 
additional contributions involving the expression $\mathcal P\partial_x^2 \phi$, exactly as in the previous subsection. 
But these contributions can be estimated perturbatively using the 
decay of the symbol for $\mathcal P$ at high frequencies.

\section{The pointwise decay properties of the nonlinear ILW equation}

In this section we consider the ILW equation with small,
localized initial data,
\begin{equation}\label{small-data}
\|\phi_0\|_{L^2} + \|x \phi_0\|_{L^2} \leq \epsilon,
\end{equation}
and seek to prove Theorem~\ref{t:quartic}, which asserts that the solution has dispersive decay up to cubic time, $|t| \lesssim \epsilon^{-2}$,
 \begin{equation}\label{want-disp}
\begin{aligned}
|\phi(t,x)| & \ \lesssim \epsilon \omega_0(t,x), 
\\
|\Tau \phi(t,x)| & \ \lesssim \epsilon \omega_1(t,x), 
\end{aligned}
\qquad t \in [0,T], \quad T \ll \epsilon^{-2}.
\end{equation}

To prove this bound in a suitable time interval $[0,T]$, it will be convenient to make a slightly weaker
bootstrap assumption, namely 
\begin{equation}\label{boot}
\begin{aligned}
|\phi(t,x)| & \ \leq C \epsilon \omega_0(t,x), 
\\
|\Tau \phi(t,x)| & \ \leq C \epsilon \omega_1(t,x), 
\end{aligned}
\qquad t \in [0,T].
\end{equation}
Here $C$ is a universal constant, which is fixed, large enough,
independent of $\epsilon$. In the proof $C$  will  not affect the 
constants in \eqref{want-disp}, but rather the choice of $T$.
Once \eqref{want-disp} is proved under the bootstrap assumption \eqref{boot},
the bootstrap assumption can be eliminated via a standard continuity argument. Throughout this section we will assume that \eqref{boot} holds.

In view of the vector field bound in Theorem~\ref{t:KS-Besov},
it would suffice to show that in the same time range we have the energy estimates
\begin{equation}
\|\phi(t)\|_{\cB^{-\frac12}_{2,\infty}} + \|L \phi(t)\|_{\cH^\frac12} \lesssim \epsilon .
\end{equation}
These two bounds are considered separately in the following two subsections.

\bigskip

\subsection{Energy bounds for \texorpdfstring{$\phi$ }{}}

Our goal here is to prove the estimate 
\begin{equation}\label{besov-u}
\|\phi(t)\|_{\cB^{-\frac12}_{2,\infty}}  \lesssim \epsilon, \qquad t \ll \epsilon^{-2}
\end{equation}
for solutions $u$ to ILW which satisfy the initial data smallness \eqref{small-data} as well as the bootstrap assumption \eqref{boot}.

\begin{remark} Based on the complete integrability of ILW and on the corresponding bounds for KdV in \cite{KT}, one would expect 
this to hold globally in time without the bootstrap assumption. However, for our purposes here we do not need such a global result, and instead we contend ourselves with a simpler, shorter argument that uses the bootstrap assumption.
\end{remark}

The $L^2$ bound for $u$ follows directly from the energy conservation, and this takes care of the high frequencies.  It remains to 
prove the low frequency bound.

We phrase the proof as a bootstrap argument, where we assume that we have the bound 
\begin{equation}\label{besov-u-boot}
\|\phi(t)\|_{\cB^{-\frac12}_{2,\infty}}  \leq C \epsilon, \qquad t \ll \epsilon^{-2},
\end{equation}
with a large constant $C$. 

Here we need to reinterpret \eqref{boot} in terms of  dyadic pointwise bounds: 
\begin{lemma}
If \eqref{boot} holds then we have
\begin{equation}
   \label{boot-d} 
   \Vert \phi_{ k} \Vert_{L^{\infty}} \leq \epsilon Ct^{-\frac{1}{2}}2^{-\frac{k}{2}} \qquad \mbox{ for } k<0,
\end{equation}
and
\begin{equation}
   \label{boot-d2} 
   \Vert \phi_{>0} \Vert_{L^{\infty}} \leq \epsilon Ct^{-\frac{1}{2}}.
\end{equation}
\end{lemma}
\begin{proof}
We begin with the bootstrap bounds \eqref{boot} and we rewrite them for the dyadic pieces $u_k$ of the solution,
\[
\vert \phi_k\vert \leq C\epsilon \omega_0 \mbox{ and } \vert \phi_k\vert \leq C\epsilon 2^{-k}\omega_1.
\]
This implies that 
\[
\vert \phi_k\vert \leq \epsilon C \min\left\{ \omega_0, 2^{-k}\omega_1\right\} ,
\]
where in  \eqref{def-omega0} and \eqref{def-omega1} the time is considered fixed.  One has
\[
\begin{aligned}
\min\left\{ \omega_0, 2^{-k}\omega_1\right\} &\leq \left\{
\begin{aligned}
& \min \left\{t^{-\frac{1}{4}}x^{-\frac{1}{4}}, t^{-\frac{3}{4}}x^{\frac{1}{4}}2^{-k}  \right\}  &\mbox{ for } x\geq -t\\
&\quad  t^{-\frac{1}{2}} &\mbox{ for } x<-t
\end{aligned}
\right. \\
&  \leq \left\{
\begin{aligned}
& t^{-\frac{1}{2}}2^{-\frac{k}{2}}
   &\mbox{ for } x\geq -t\\
&\quad  t^{-\frac{1}{2}} &\mbox{ for } x<-t.
\end{aligned}
\right. \\
\end{aligned}
\]
Hence, for $k<0$, we obtain
\[
\vert  \phi_k \vert \leq \epsilon C t^{-\frac{1}{2}}2^{-\frac{k}{2}},
\]
which implies the bound in \eqref{boot-d}. The bound in \eqref{boot-d2} also follows.
\end{proof}

Given a frequency $k < 0$, we need to show that  
\begin{equation}\label{want-k}
\| \phi_k\|_{L^2} \lesssim 2^{\frac{k}2} \epsilon.
\end{equation}
First we verify that this holds at the initial time,
\[
\| \phi_k(0)\|_{L^2} \lesssim 2^{\frac{k}2}(\|\phi(0)\|_{L^2} + \|x\phi(0)\|_{L^2}).
\]
Then we want to propagate this in time. 
For this we compute
\begin{equation}
\label{boot-propagate}
\frac{d}{dt} \|\phi_k\|_{L^2}^2 =  \int  \phi_k \cdot \partial_x P_k (\phi^2)\, dx  .
\end{equation}
The $\partial_x$  gives the desired $2^k$ factor,
so we can use our bootstrap assumption to estimate the integral on the right by
\[
\leq \epsilon C 2^{\frac{3k}2} \|P_k (\phi^2)\|_{L^2}.
\]
It remains to bound in $L^2$ the expression $P_k(\phi^2)$.  Here we use the Littlewood-Paley decomposition 
\[
P_k(\phi^2)\approx \phi_{\leq k}\phi_k +\sum_{k \leq j < 0} P_k(\phi_j\phi_j) +P_k(\phi_{>0}\phi_{>0}),
\]
and the bootstrap bounds for $\phi$
in both $L^2$ (see \eqref{besov-u-boot}) and $L^\infty$ (see \eqref{boot-d}) as follows
\[
\begin{aligned}
\Vert P_k(\phi^2)\Vert_{L^2}&\leq \Vert \phi_{\leq k}\Vert_{L^2} \Vert \phi_k\Vert_{L^{\infty}}+\sum_{k\leq j< 0} \Vert \phi_{j}\Vert_{L^2} \Vert \phi_{j}\Vert_{L^{\infty}} + \Vert \phi_{>0}\Vert_{L^2} \Vert \phi_{> 0}\Vert_{L^{\infty}} \\
&\leq 2^{\frac{k}{2}} \Vert \phi_k\Vert_{\cB^{-\frac12}_{2,\infty}}\epsilon Ct^{-\frac{1}{2}}2^{-\frac{k}{2}}+ \sum_{k\leq j< 0} 2^{\frac{j}{2}}\Vert \phi_j\Vert_{\cB^{-\frac12}_{2,\infty}}\epsilon Ct^{-\frac{1}{2}} 2^{-\frac{j}{2}}  + \Vert \phi_{>0}\Vert_{\cB^{-\frac12}_{2,\infty}} \epsilon Ct^{-\frac{1}{2}}\\
&\leq (\epsilon Ct^{-\frac{1}{2}}  + \epsilon Ct^{-\frac{1}{2}}k)
\Vert \phi\Vert_{\cB^{-\frac12}_{2,\infty}}
\\
& \leq \epsilon^2 C^2 t^{-\frac{1}{2}} (2+k).
\end{aligned}
\]
Returning to  the bound in \eqref{boot-propagate}, it follows that
\[
\frac{d}{dt} \| \phi_k\|_{L^2}^2 \leq C\epsilon^32^{\frac{3k}{2}}t^{-\frac{1}{2}} k.
\]
Integrating this on the time interval $[0,T]$ shows that  the desired bound \eqref{want-k}  propagates in time for as long as $t\ll \epsilon^{-2}$.

\subsection{Energy bounds for \texorpdfstring{$Lu$}{}}
Our goal here is to prove the estimate 
\begin{equation}\label{l-ee}
\|L \phi(t)\|_{\cH^\frac12} \lesssim \epsilon, \qquad t \ll \epsilon^{-2},
\end{equation}
for solutions $\phi$ to ILW which satisfy the initial data smallness \eqref{small-data} as well as the bootstrap assumption \eqref{boot}.

The difficulty is that one cannot obtain directly 
a propagation bound for $L\phi$. While $L$ commutes with the linear flow, its action on the nonlinear flow is more 
difficult to track; in particular, the operator $L$
is not naturally associated to any symmetry of the nonlinear flow,
so one cannot directly write a good evolution equation for $L\phi$.

To address the above mentioned difficulty, the idea
is to find an operator $L^{NL} $, which is a nonlinear
correction of $L$, and so that $L^{NL} \phi$ satisfies 
good energy estimates.  This idea is analogous to the one found in \cite{IT-BO}, though the analysis here
is more difficult due in particular to the lack of scale invariance for the ILW flow. Once this correction is found so that we have the 
good energy estimates
\begin{equation}\label{LNL-ee}
\|L^{NL} \phi(t)\|_{\cH^\frac12} \lesssim \epsilon,
\end{equation}
we will return to \eqref{l-ee} by directly estimating the correction.

 The operator $L^{NL}$ is obtained by adding  a quadratic correction to the linear operator $L$,
 \begin{equation}\label{def:LNL}
 L^{NL} \phi = L \phi+ tB(\phi,\phi),
 \end{equation}
 where $B$ is a symmetric translation invariant bilinear form.
Defining the new variable 
\begin{equation}\label{def-v}
 v = L\phi + tB(\phi,\phi),   
\end{equation}
it is easily seen that $v$ must solve an equation of the form
\begin{equation}\label{v-eqn}
Pv = C(\phi,v) + t R(\phi,\phi,\phi) + D(\phi,\phi),
\end{equation}
where $C$, $R$ and $D$ are also translation invariant multilinear forms.

Ideally, one would like to find a good choice of $B$, so that $C$ is antisymmetric  as a linear operator
acting on $v$, and so that $R=0$, $D=0$. As proved in \cite{IT-BO}, this is indeed the case for Benjamin-Ono, which approximates ILW at high frequency; we will review this below. However, this is not the case  for KdV.  So we will have to compromise somewhat for ILW. In particular, instead of asking for $C$ to be antisymmetric, we will ask to have good symbol bounds for its symmetric part
\[
C_{sym} = \frac12 (C+ C^*).
\]
To understand what is the symmetric part of $C$ as a linear operator acting on $v$ we compute
\[
0=\int C(u,v) w + C(u,w) v \, dx = \int_{\xi+\eta+\zeta = 0}
(c(\xi,\eta) + c(\xi,\zeta)) \, \hat u(\xi) \hat v(\eta) \hat w(\zeta) \, dA,
\]
which leads to the symbol expression 
\begin{equation}\label{C-antisym}
c_{sym}(\xi,\eta) = \frac12(c(\xi,\eta) + c(\xi,\zeta))  \qquad \xi+\eta+\zeta = 0.
\end{equation}
As we will show next, a good  correction $B$ can be computed, but the analysis is not entirely straightforward.  Our choice of $B$ is described in the following proposition:

\begin{proposition}\label{p:B}
There exists a correction $B$ for the operator $L$ associated to ILW so that the following properties hold for $B$, and $C,R,D$ in \eqref{v-eqn}:
\begin{enumerate}[label=(\roman*)]
    \item Bounded $B$: the symbol for $B$ is smooth, 
    real, symmetric, even, with symbol type regularity, and satisfies
   \begin{equation}\label{b-size}
  |b(\xi,\eta)| \lesssim 1.     
   \end{equation} 
  \item Almost antisymmetric $C$: the symbol for $C$ is smooth, 
  purely imaginary, odd, with symbol type regularity, and 
  satisfies
  \begin{equation}\label{c-size}
  |c(\xi,\eta)| \lesssim |\xi|+|\eta|   .  
   \end{equation} 
 In addition, the  symbol for its antisymmetric part satisfies
   \begin{equation}\label{csym-size}
  |c_{sym}(\xi,\eta)| \lesssim |\xi| |\eta| e^{-c(|\xi|+|\eta|)}.      
   \end{equation} 
\item Bounded $R$, $D$: the symbols for $R$ and $D$ are smooth and satisfy
\begin{equation}\label{r-size}
  |r(\xi,\eta,\zeta)| \lesssim |\Tau(\xi)| + |\Tau(\eta)| + |\Tau(\zeta)|,     
   \end{equation} 
  \begin{equation}\label{d-size}
  |d(\xi,\eta)| \lesssim 1,     
   \end{equation}  
both with symbol type regularity.   
\end{enumerate}
\end{proposition}
Here we need to be careful what we call symbol type regularity. In the case of $B$, for instance, we have also some nontrivial dependence on $\xi+\eta$, and similar considerations apply for $C$, $R$ and $D$. This will be described in greater detail in the Lemmas which we establish within the proof of the Proposition.

We also remark on the parity conditions for these symbols, which are set to guarantee that for real valued $u$,
all the outputs are also real valued functions.

To simultaneously find the symbols of $B$ and $C$,
we perform an algebraic computation that sets the quadratic terms in the $v$ equation \eqref{v-eqn} to zero. Thus
\begin{equation}
\label{Lnl}
\begin{aligned}
Pv=&P(L\phi+tB(\phi,\phi))\\
=&P(L\phi)+tP(B(\phi,\phi)) +B(\phi,\phi)\\
=&L(P\phi)+tP(B(\phi,\phi)) +B(\phi,\phi)\\
=&L(\phi\phi_x)+t\left[ 2B(\phi_t, \phi)+A(D)B(\phi,\phi)\right] + B(\phi,\phi)\\
=&L(\phi\phi_x)+t\left[ - 2B(A(D)\phi, \phi)+A(D)B(\phi,\phi) \right] + B(\phi,\phi)+ 2tB(\phi,\phi\phi_x).\\
\end{aligned}
\end{equation}

Comparing the quadratic terms here with those in \eqref{v-eqn} we arrive at the quadratic identity
\begin{equation}\label{quad-id}
 C(\phi,L\phi) + D(\phi,\phi) = L(\phi\phi_x)+t\left[ - 2B(A(D)\phi, \phi)+A(D)B(\phi,\phi) \right] + B(\phi,\phi).
\end{equation}
If this holds, then for $v$ we indeed obtain the equation 
\eqref{v-eqn}, with 
\begin{equation}\label{def-R}
R(\phi,\phi,\phi) = 2 B(\phi,\phi\phi_x) - C(\phi, B(\phi,\phi)).
\end{equation}

Now our task is to find the translation invariant bilinear forms $B$ and $C$,
described by their symbols $b(\xi,\eta)$ and $c(\xi,\eta)$,
so that the relation \eqref{quad-id} holds.

If we identify the $x$ terms on the left and on the right 
in \eqref{quad-id} then we get the leading order relation
\begin{equation}\label{first-C}
C(\phi,\phi) = \phi \phi_x,
\end{equation}
as well as the secondary relation 
\begin{equation}\label{first-C+}
D(\phi,\phi) = B(\phi,\phi) - i C_\eta(\phi,\phi) ,    \end{equation}
arising from commuting $x$ outside $C$. 

The first relation provides us with the symmetric part of $C$. Precisely,  if we split the bilinear form $C$
into a symmetric and an antisymmetric\footnote{Note that this is not the same as having $C$ antisymmetric as an operator acting on $v$.} part,
\[
C = C^s + C^a,
\]
then  the first equation \eqref{first-C}
gives
\[
C^s(\phi,v) = \frac12 \partial_x(\phi v) ,
\]
or at the symbol level,
\begin{equation}
\label{e:c-sym}
c^s(\xi,\eta) = \frac{i}{2}(\xi+\eta).
\end{equation}
The secondary relation \eqref{first-C+}, on the other hand,
may be thought of as the definition of $D$, once
$C$ is determined; here we remark that $D$ should be thought of as a symmetric bilinear form, so the symbol $C_\eta$
needs to be symmetrized.

 Next we identify the $t$ terms in \eqref{quad-id} in order to
get a second relation, 
\begin{equation}\label{second-C}
i C(\phi,A_\xi(D) \phi) = i A_\xi(D)(\phi \phi_x) + \left[ - 2B(A(D)\phi, \phi)+A(D)B(\phi,\phi) \right]. 
\end{equation}
We substitute here the symmetric part of $C$, computed above, to obtain
\[
 i C^a(\phi,A_\xi(D) \phi) + \frac{i}2 \partial_x(\phi A_\xi(D)\phi) =  \frac{i}2 \partial_x A_\xi(D)(\phi^2) + \left[ - 2B(A(D)\phi, \phi)+A(D)B(\phi,\phi) \right]. 
\]
Now we write the corresponding relation for the symbols, making sure we symmetrize on the left.
This gives
\[
\frac{i}2 c^a(\xi,\eta) ( a_\xi(\eta) - a_\xi(\xi)) - \frac{1}{4}(\xi+\eta)(a_\xi(\xi) + a_\xi(\eta))
\!= \!b(\xi,\eta) (a(\xi+\eta)-a(\xi)-a(\eta)) -\frac12 (\xi+\eta) a_\xi(\xi+\eta).
\]
We want to think of this as an equation for $b$, so we rewrite it as
\[
\begin{aligned}
b(\xi,\eta) (a(\xi+\eta)-a(\xi)-a(\eta)) &= \frac12 \left[ (\xi+\eta) a_\xi(\xi+\eta) - \xi a_\xi (\xi) - \eta a_\xi(\eta)\right]
\\
&\qquad + \left( \frac{i}2 c^a(\xi,\eta) - \frac14(\xi- \eta)\right) ( a_\xi(\eta) - a_\xi(\xi)) .
\end{aligned}
\]
We remark that so far, this does not uniquely determine
$c^a(\xi, \eta)$ and $b (\xi, \eta)$; instead, it allows us to compute $b(\xi, \eta)$
given $c^a(\xi, \eta)$. Observe also that the coefficient of $b(\xi, \eta)$ corresponds to the quadratic resonance relation that we know only vanishes when one of the entry frequencies is zero ($\xi$ or $\eta$), or when the outcome frequency $\xi + \eta$ is zero. The smoothness of the symbol $b(\xi,\eta)$ on this resonance set is very important, so it is essential to choose $c^a(\xi, \eta)$ that the right hand side above vanishes on the set of quadratic resonances. Here we recall that 
\[
a(\xi)=i\xi^2\coth(\xi).  
\]
To simplify the computations we separate $b$ into two parts 
\[
b(\xi, \eta) := \frac12 b_1 (\xi, \eta)+b_2 (\xi, \eta),
\]
where $b_1$ is given by
\begin{equation}\label{def-B1}
{b_1}(\xi, \eta):= \frac{a'(\xi +\eta) (\xi +\eta)- a'(\xi)\xi-a'(\eta)\eta}{a(\xi +\eta)-a(\xi) - a(\eta)}
\end{equation}
is easily seen to be an everywhere smooth zero order 
symbol.  Here the derivative $'$ denotes the derivative with respect to $\xi$. We remark on the asymptotic limits for $b_1$,
\begin{equation}
b_1 \approx \left\{ 
\begin{array}{ll}
2 & \qquad |\xi|+|\eta| \to \infty \cr 
3 &  \qquad |\xi|+|\eta| \to 0 .
\end{array}
\right.
\end{equation}
These correspond to the Benjamin-Ono, respectively the KdV regimes.

At this point we are left with an equation connecting $b_2$ and $c^a$, namely
\begin{equation}
\label{b2-def}
b_2(\xi,\eta) (a(\xi+\eta)-a(\xi)-a(\eta)) =  \left(\frac{i}2 c^a(\xi,\eta)- \frac14(\xi-\eta)\right) ( a_\xi(\eta) - a_\xi(\xi)).
\end{equation}
To respect the parity symmetries in our problem, 
here we assume that $c^a$ is purely imaginary and 
odd, so that $b_2$ is real, symmetric and even.
In order to have a smooth division, we must insure that $c^a$ is chosen so that 
the real expression
\begin{equation}\label{tca}
  \tilde c^a(\xi,\eta) := \frac{i}2 c^a(\xi,\eta)- \frac14(\xi-\eta)
\end{equation}
vanishes at $\xi = 0$ and at $\eta = 0$.
Still, this does not uniquely determine
the choice of $\tilde c^a$, so we need to 
make a careful selection.

Our choice should match the KdV choice at low frequency, and the Benjamin-Ono choice at high frequency. Because of this, it is useful
to stop and examine these two choices:
\bigskip 

\textbf{ The KdV case:}
Here, following the KdV set-up in \cite{KIT},  the choice would be to take 
\begin{equation}\label{kdv}
b_2 (\xi, \eta)=0,\qquad \tilde c^a(\xi,\eta) = 0, \qquad c^a(\xi,\eta) = -\frac{i}{2}(\xi-\eta).
\end{equation}
While this would be reasonable at low frequency,  it would not be
satisfactory in general. This is because, ideally, the other requirement is that $C$ is almost antisymmetric as a linear operator; otherwise, the $C$ term yields  a nontrivial contribution in the energy estimates for $v$, of the form  $\int \phi_x v^2\, dx$, which is difficult to control and we want to avoid in the high frequency limit.

\bigskip

\textbf{ The Benjamin-Ono case:} This corresponds to the analysis in \cite{IT-BO}. Here $\tilde c^a_{IT-BO}$ should be homogeneous of order one. It suffices to confine our choice to piecewise linear $\tilde c^a_{IT-BO}$, and then the requirement that it vanishes on the axes 
$\xi=0$ and $\eta=0$ and on the diagonal $\xi=\eta$
lead to a unique choice in the first and the third quadrant,
\begin{equation}\label{first-quad}
\tilde c^a_{IT-BO}(\xi,\eta) = 0,  \quad \mbox{when }\quad \sgn \xi = \sgn \eta .
\end{equation}
To determine $\tilde c^a_{IT-BO}$ in the remaining quadrants 
we additionally impose the requirement that $C$ be antisymmetric, 
$C_{sym} = 0$. We set the symbol $c_{sym}$ in \eqref{C-antisym}
to zero, and rewrite it as a condition for the antisymmetric part $c^a$ of $c$. We get
\[
c^a_{IT-BO}(\xi,\eta) + c^a_{IT-BO}(\xi,\zeta) =  \frac{i}{2} (\eta+\zeta) = -\frac{i\xi}2,
\]
and further 
\begin{equation}\label{C-sym}
\tilde c^a_{IT-BO}(\xi,\eta) + \tilde c^a_{IT-BO}(\xi,\zeta) = 
 \frac{\xi}4 - \frac{3\xi}{4} = -\frac{\xi}2.
\end{equation}
Combining this with \eqref{first-quad} gives 
\[
\tilde c^a_{IT-BO}(\xi,-\xi-\eta) = - \frac{\xi}2, \quad  \mbox{ when }\quad \sgn \xi = \sgn \eta,
\]

or equivalently, using also the antisymmetry condition \eqref{C-sym},
\[
\tilde c^a_{IT-BO} (\xi,\eta) =- \frac{\xi-\eta}4 + \frac{\xi+\eta}{4} \sgn(\xi+\eta) \sgn \xi,
\qquad \sgn \xi \neq \sgn \eta.
\]
For later use we further write it in the equivalent form
\begin{equation} \label{tca-bo}
\begin{aligned}
  \tilde c^a_{IT-BO} (\xi,\eta) = & \ 
\frac12(1-\sgn \xi \sgn\eta) \left( - \frac{\xi-\eta}4 + \frac{\xi+\eta}{4} \sgn(\xi+\eta) \sgn \xi
\right)
\\
= &- \frac14  \xi \sgn \eta(\sgn(\xi+\eta) - \sgn \xi) +\frac14 \eta \sgn \xi (\sgn(\xi+\eta) - \sgn \eta).
\end{aligned}  
\end{equation}

Hence in the Benjamin-Ono equation's case we get
\[
\begin{aligned}
-i c_{IT-BO}(\xi,\eta) = & \ \eta -2   \tilde c^a_{IT-BO}(\xi,\eta)
\\
= & \ \eta  + \left(1 - \sgn \xi \sgn \eta\right)  \left( \frac{\xi-\eta}4 - \frac{\xi+\eta}{4} \sgn(\xi+\eta ) \sgn \xi \right)
\\
= & \ \frac14 (\xi+3\eta) + \frac14 (\eta-\xi)\sgn \xi \sgn \eta  - \frac14 (\xi+\eta) \sgn \xi \sgn(\xi+\eta) + \frac14(\xi+\eta) \sgn(\xi+\eta) \sgn \eta
\\
= & \ \frac14 (\xi+2\eta) -\frac14 \xi \sgn \xi \sgn \eta  - \frac14 \xi \sgn \xi \sgn(\xi+\eta) + \frac14(\xi+2\eta) \sgn(\xi+\eta) \sgn \eta .
\end{aligned}
\]
Examining the last expression, we see that the operator $C_{IT-BO}$ is antisymmetric and can be written in the form 
\begin{equation}\label{C-BO}
4 C_{IT-BO}(\phi, v) = (\phi \partial_x + \partial_x \phi)v - H  (\phi \partial_x + \partial_x \phi)Hv + (H \partial_x \phi)H v
+ H [(H \partial_x \phi) v].
\end{equation}
We also compute the corresponding correction $B$ (we will do so at the level of the symbol $b$).
As noted earlier, we have $b_{BO,1}=2$, while for $b_{BO,2}$ we obtain 
\[
b_{BO,2}(\xi,\eta) = -\frac14 (1 - \sgn \xi \sgn \eta)  \frac{(\xi-\eta - (\xi+\eta) \sgn(\xi+\eta) \sgn \xi )(|\eta|-|\xi|)}{ (\xi+\eta)|\xi+\eta| - \xi |\xi| - \eta|\eta|}.
\]
The first factor restricts its support to the region where $\xi$ and $\eta$ have opposite signs, so we 
can further rewrite it as
\[
\begin{aligned}
b_{BO,2}(\xi,\eta) = & \ -\frac14 (1 - \sgn \xi \sgn \eta)  \frac{(\xi-\eta - (\xi+\eta) \sgn(\xi+\eta) \sgn \xi )\sgn \eta(\xi+\eta)}{ (\xi+\eta)|\xi+\eta| - \sgn \xi(\xi^2-\eta^2)} 
\\ 
= & \
-\frac14 (1 - \sgn \xi \sgn \eta)  \frac{(\xi-\eta)\sgn \eta + (\xi+\eta) \sgn(\xi+\eta)}{ |\xi+\eta| - \sgn \xi(\xi-\eta)} 
\\ 
= & \
-\frac14 (1 - \sgn \xi \sgn \eta) .
\end{aligned}
\]
Hence we obtain
\[
b_{IT-BO}(\xi,\eta) = \frac34 + \frac14 \sgn \xi \sgn \eta.
\]
In operator form, this reads
\[
B_{IT-BO}(\phi,\phi) = \frac14 ( 3 \phi^2 - (H\phi)^2).
\]
We carefully observe that this is exactly the correction 
used in \cite{IT-BO}. Finally, we note that 
a direct computation shows that $R=0$ and $D=0$ in the Benjamin-Ono equation case.

\bigskip

Now we return to the choice of $\tilde c^a$ in the 
ILW case. We seek to make a smooth choice, but which matches the Benjamin-Ono choice in the high frequency limit, and the KdV choice at  low frequency. Inspired by the second expression in \eqref{tca-bo},
we will set 
\begin{equation} \label{tca-ilw}
\begin{aligned}
  \tilde c^a (\xi,\eta) = & \ \frac14  \xi \Tau (\eta)(\Tau(\xi+\eta) - \Tau(\xi)) -\frac14 \eta \Tau (\xi) (\Tau(\xi+\eta) - \Tau (\eta)).
\end{aligned}  
\end{equation}
As needed, this is smooth, real, odd and antisymmetric,
and vanishes on the axes $\xi=0$ and $\eta=0$. 
Hence the division in \eqref{b2-def} yields a smooth,
real, even and symmetric symbol $b_2$.
Further, we note that $b_2(0,0)= 0$, which agrees 
with the KdV choice. This symbol clearly satisfies the 
bound \eqref{c-size} and the corresponding regularity 
requirements.

Next, we check whether $C$ is antisymmetric. For this we compute, with $\xi+\eta+\zeta = 0$, 
\[
\begin{aligned}
4 (\tilde c^a(\xi,\eta) + \tilde c^a(\xi,\zeta))
= & \  \xi \Tau (\eta)(\Tau(\zeta) + \Tau(\xi)) -\eta \Tau (\xi) (\Tau(\zeta) + \Tau (\eta))
\\
& \ +\xi \Tau (\zeta)(\Tau(\eta) + \Tau(\xi)) -\zeta \Tau (\xi) (\Tau(\eta) + \Tau (\zeta))
\\
= & \  2 \xi ( \Tau(\xi) \Tau(\eta) + \Tau(\xi) \Tau(\zeta)
+\Tau(\eta) \Tau(\zeta))
\\
= & \ 2\xi + 2 \xi (\Tau(\xi) \Tau(\eta) + \Tau(\xi) \Tau(\zeta)
+ \Tau(\eta) \Tau(\zeta)-1).
\end{aligned}
\]
Thus, returning to $c$, we get
\[
c(\xi,\eta) + c(\xi,\zeta) = \frac{i}{2} \xi (1-\Tau(\xi) \Tau(\eta) - \Tau(\xi) \Tau(\zeta)
- \Tau(\eta) \Tau(\zeta)).
\]
Here we do not get full cancellation in the last bracket, but  we do get rapid decay at high frequency. To describe this decay, it is convenient to view $c$ as a function on the plane $\xi+\eta+\zeta=0$. Using the notation
\[
 \xi_{hi} := \max\{|\xi|,|\eta|,|\zeta|\} ,
 \qquad  \xi_{lo} := \min\{|\xi|,|\eta|,|\zeta|\} ,
\]
we summarize the information about $c$ in the next lemma:

\begin{lemma} \label{l:C-sym}
The symbol $c$ is given by 
\begin{equation}
c(\xi,\eta) = i \eta - 2i \tilde c^a(\xi,\eta)    
\end{equation}
with $\tilde c^a$ defined by \eqref{tca-ilw}. The symbol of its symmetric part of $C$ has the form
\begin{equation}
c_{sym}(\xi,\eta) = c(\xi,\eta)+c(\xi,\zeta)    
= \xi g(\xi,\eta),
\end{equation}
where the symbol $g$ is smooth and exponentially decaying, 
\[
|g(\xi,\eta)| \lesssim e^{-c \xi_{hi}}
\]
along with its derivatives.
\end{lemma}

We next consider the symbol $b$. By definition the symbol $\tilde c^a$ 
vanishes at $\xi = 0$ and $\eta = 0$, which guarantees  a smooth division in \eqref{b2-def}. It remains to consider
the high frequency asymptotics for $b$, which  are described next:
\begin{lemma}\label{l:B}
The symbol $b(\cdot, \cdot)$ has the representation
\begin{equation}
b(\xi,\eta) = b^\sharp(\xi,\eta) + b^r(\xi,\eta)  ,  
\end{equation}
where the leading part $b^\sharp$ is explicit,
\begin{equation}
b^\sharp(\xi,\eta) = \frac34 - \frac14 \Tau(\xi) \Tau(\eta) ,
\end{equation}
and the residual part decays,
\begin{equation}\label{br}
|b^r(\xi,\eta)| \lesssim \langle \xi_{hi} \rangle^{-1} e^{-c \xi_{lo}}.
\end{equation}
along with its derivatives.
\end{lemma}
We comment further on the structure of $b^r$, which is 
not in a classical symbol class. Again, this is best seen 
as a function of $(\xi,\eta,\zeta)$ on the plane 
$\xi+\eta+\zeta = 0$, and decays exponentially away from the lines $\xi = 0$, $\eta = 0$, respectively $\zeta = 0$.
Near each of these axes we obtain an asymptotic description of $b^r$. Near $\eta = 0$ for instance,
we have an asymptotic expansion of the form
\begin{equation}\label{expansion}
b^r(\xi,\eta) \approx \sum_{j \geq 1} 
b^{j1}(\eta) \xi^{-j} + b^{j2}(\eta) \xi^{-j} \sgn \xi
\end{equation}
with exponentially decaying coefficients 
$b^{j1}$, $b^{j2}$. This description suffices in order 
to obtain good kernel bounds on the associated  operators.

\begin{proof}
We consider separately the components $b_1$ and $b_2$ given by \eqref{def-B1}, respectively \eqref{b2-def}.
For $b_1$ we have its Benjamin-Ono counterpart $b_{BO,1} = 2$. Subtracting it, we write 
\[
b_1(\xi,\eta)   = 2+ 
\frac{m(\zeta) + m(\xi)+m(\eta)}{a(\zeta)+a(\xi) + a(\eta)}, \qquad \xi+\eta+\zeta = 0.
\]
where the auxiliary function
\[
m(\xi) := \xi a'(\xi) - 2a(\xi)
\]
is odd and exponentially decaying at infinity.
This guarantees that the division is smooth, so we consider its asymptotics at infinity:

\begin{enumerate} [label=\alph*)]
    \item  If $\xi,\eta,\zeta$ are all large and comparable then we get exponential decay from the numerator.
 \item Else, by symmetry we consider the case when $|\zeta| \lesssim |\xi| \approx |\eta|$. Then, modulo exponentially decaying tails we can write
\[
b_1 - 2 \approx \frac{m(\zeta)}{\xi|\xi| + \eta|\eta| + a(\zeta)} = \frac{m(\zeta)}{\zeta \sgn \xi (2\xi+\zeta + \zeta^{-1} a(\zeta) \sgn \xi )}.
\]
Here $\zeta^{-1} m(\zeta)$ is smooth, even, vanishing at zero and exponentially decaying at infinity.
On the other hand the remaining expression  is a classical symbol of order $-1$,
which for  $|\zeta| \ll |\xi|$ we can expand in a Taylor series
\[
\frac{1}{2\xi+\zeta + \zeta^{-1} a(\zeta) \sgn \xi} = 
\frac{1}{2\xi} - \frac{\zeta + \zeta^{-1} a(\zeta) \sgn \xi }{4 \xi^2} + \cdots
\]
The powers of $\zeta$ are absorbed by the $\zeta^{-1} m(\zeta)$ factor, so this expansion yields 
better and better errors for the product. Here one can do a similar expansion with $\eta$ and $\xi$ interchanged.
\end{enumerate}
The above discussion shows that the contribution of $b_1-2$ is entirely placed into $b^r$.

\medskip

Now we turn our attention to $b_2$, for which we have 
\[
b_2 =  \frac{ \tilde c^a(\xi,\eta) ( a'(\xi) - a'(\eta))}{a(\zeta)+a(\xi)+a(\eta)}.
\]
Again by the choice of $\tilde c^a$ we have a smooth division, so it suffices to consider the high frequency asymptotics.  Due to a lack of full symmetry, we now consider three cases. 

\begin{enumerate} [label=\alph*)]
 \item  If $\xi,\eta,\zeta$ are all large and comparable then we directly get exponential decay to the Benjamin-Ono  setting.

\item If  $\xi$ and $\eta$ are the high frequencies,
i.e. $|\zeta| \ll |\xi|,|\eta|$, then $\xi$ and $\eta$ 
have opposite signs. Modulo exponentially decaying factors we write
\[
 \tilde c^a (\xi,\eta) \approx -\frac14 (\xi - \eta)
 - \frac{1}4 (|\xi| - |\eta|) \tanh{\zeta} 
 \]
Then for $b_2$ we arrive at
\[
\begin{aligned}
b_2(\xi,\eta) \approx & \ \frac{[ (\xi - \eta)
 -  (|\xi| - |\eta|) \tanh \zeta]
 (|\xi|-|\eta|)}{2 \zeta ( |\xi - \eta| + \zeta^{-1} a(\zeta))} 
\\
 \approx & \ \frac{ (\xi - \eta)
 + \zeta \tanh \zeta \sgn \xi
 }{2  ((\xi - \eta)+ \zeta^{-1}a (\zeta) \sgn \xi) } 
\\
\approx & \ \frac12\left(1 + \frac{m_1(\zeta)\sgn \xi}{ 2\xi+\zeta +  \zeta^{-1} a(\zeta) \sgn \xi}\right) 
\end{aligned}
\]
where 
\[
m_1(\zeta) = \zeta\tanh{\zeta} - \zeta^{-1} a(\zeta)
\]
decays exponentially at infinity. The $1/2$ term agrees
with the corresponding Benjamin-Ono contribution, while the last term can be placed into $b^r$ exactly as in the $b_1$ case, by taking a Taylor expansion.

\item The last case is when  $\xi$ and $\zeta$ are the high frequencies, i.e. $|\eta| \ll |\xi|, |\zeta|$. 
By symmetry the same applies when $\eta$ and $\zeta$ 
are the small frequencies. Here we write, again with exponentially small errors vanishing at $\eta = 0$,
\[
\tilde c^a(\xi,\eta) \approx \frac14 \eta (1-\sgn \xi \tanh{\eta}).
\]
This allows us to compute with exponentially small errors
\[
\begin{aligned}
b_2(\xi,\eta) \approx  & \ \frac{ \tilde c^a(\xi,\eta) ( 2|\xi| - a'(\eta))}{-\eta|\xi-\zeta| + a(\eta)}
\\
\approx  & \ -\frac14 \frac{ (1 - \sgn \xi \tanh{\eta})(2|\xi| - a'(\eta))}{ |\xi - \zeta| - \eta^{-1} a(\eta)} \\ 
 \approx  & \ - \frac14 \frac{ (1 - \sgn \xi \tanh{\eta})(2\xi  - \sgn \xi a'(\eta))}{ 2\xi + \eta - \eta^{-1} a(\eta) \sgn \xi}
\\
\approx & \ - \frac14 (1-\sgn{\xi} \tanh{\eta})\left[ 1- 
\frac{  \eta - (\eta^{-1} a(\eta) - a'(\eta)) \sgn \xi
}{ 2\xi + \eta + \eta^{-1} a(\eta) \sgn \xi }\right]\, .
\end{aligned}
\]
Here the contribution of the first term is the leading part. To place the second term into the residual term $b^r$ we need to 
check that it has exponential decay in $\eta$ when $|\eta|\ll |\xi|$. This comes from the factor
$1-\sgn{\xi} \tanh{\eta}$ when $\xi$ and $\eta$ have the same sign, and from the factor $ \eta - (\eta^{-1} a(\eta) - a'(\eta)) \sgn \xi$ otherwise.

\end{enumerate}

\end{proof}

Once we have a good understanding of $b$, we are able to obtain bounds for $d(\cdot, \cdot)$, which is the  symbol associated to the bilinear form $D$ given in \eqref{first-C+}:

\begin{lemma}\label{l:D}
The symbol $d$ satisfies
\begin{equation}
   |d(\xi,\eta)| \lesssim  e^{-c \xi_{lo}}   .  
\end{equation}

\end{lemma}
As in the case of the residual term $b^r$, we remark that 
at high frequency the symbol $d$ is concentrated near the axes $\xi=0$, $\eta = 0$ respectively $\zeta = 0$, and that in those regions is admits expansions similar to 
\eqref{expansion}, but starting at order $0$ rather than $-1$.

\begin{proof}
We recall that, by \eqref{first-C+}, 
\[
d(\xi,\eta) = b(\xi,\eta) - i [c_\eta(\xi,\eta)]^{s} =
b(\xi,\eta) - 1 + 2 [\tilde c^a_\eta(\xi,\eta)]^{s}.
\]
where the superscript "s" denotes the symmetrization
of the symbol. By the antisymmetry of $\tilde c^a$ we have 
\[
\begin{aligned}
2 [\tilde c^a_\eta(\xi,\eta)]^{s} = & \ (\partial_\eta-\partial_\xi) \tilde c^a(\xi,\eta)
\\
= & \  \tanh{\xi}(\tanh{(\xi+\eta)}-\tanh{\eta})
+ \tanh{\eta}(\tanh{(\xi+\eta)}-\tanh{\xi})
\\ & \ - \eta ( \tanh{\xi} \tanh'\eta + \tanh'{\xi} (\tanh{(\xi+\eta)}-\tanh{\eta}))
\\ & \ - \xi ( \tanh{\eta} \tanh'\xi + \tanh'{\eta} (\tanh{(\xi+\eta)}-\tanh{\xi})).
\end{aligned}
\]
The terms on the last two lines are easily seen to have the desired boundedness and exponential decay. 
Hence, using also Lemma~\eqref{l:B}, we write
\[
\begin{aligned}
d(\xi,\eta) = & \ -\frac14 + \frac14 \tanh \xi \tanh \eta
 + \frac12  \tanh \xi (\tanh(\xi+\eta) - \tanh \eta)
\\
& + \frac12 \tanh \eta (\tanh(\xi+\eta) - \tanh \xi)]
+ O(e^{-c \xi_{lo}})
\\
 = & \  \frac14 (-1 -\tanh{\xi} \tanh{\eta}+ \tanh{\xi} \tanh(\xi+\eta) + \tanh{\eta} \tanh(\xi+\eta)) + O(e^{-c \xi_{lo}})
\\
 = & \ O(e^{-\xi_{lo}}),
\end{aligned}
\]
as needed.
\end{proof}
Here we remark on the Benjamin-Ono counterpart of this computation, which is exact with $\tanh$ replaced by $\sgn$, without any errors.
This yields $d_{IT-BO}=0$.

Finally, we consider the trilinear form $R$, which we compare with its 
Benjamin-Ono's counterpart:

\begin{lemma}\label{l:R}
  The symbol $r$ is uniformly smooth and bounded, with  
$r(0,0,0) = 0$, and admits a representation 
\begin{equation}
r(\xi,\eta,\zeta) = \Tau(\xi) s_1(\xi,\eta,\zeta)
+ \Tau(\eta) s_2(\xi,\eta,\zeta) + \Tau(\zeta) s_3(\xi,\eta,\zeta),
\end{equation}
where the multilinear operators $S_j$ have integrable
and rapidly decreasing kernels.
\end{lemma}
Here we remark that the quotients $s_1,s_2,s_3$ in this decomposition are concentrated along the planes 
$\xi=0$, $\eta= 0$, $\xi+\eta = 0$, $\xi+\zeta = 0$,
$\eta+\zeta = 0$, with exponential decay away from these 
planes and asymptotic expansions similar to 
\eqref{expansion} along these sets.

\begin{proof}
We recall the formula for the trilinear form $R(\cdot, \cdot, \cdot)$ given in \eqref{def-R}:
\[
R(\phi,\phi,\phi) = 2 B(\phi,\phi_x) - C(\phi, B(\phi,\phi)).
\]
Here we use $c(\xi,\eta) = i(\eta - 2 \tilde c^a(\xi,\eta))$
to rewrite this as 
\[
R(\phi,\phi,\phi) = 2 B(\phi,\phi\phi_x) - 2 \phi B(\phi,\phi_x) - 2i \tilde C^a(\phi, B(\phi,\phi)).
\]
For $B$ we can use Lemma~\ref{l:B} to reduce the problem to the case when $B=B^\sharp$ as the derivatives in the first two terms as well as the one present in $C$
are cancelled by the one derivative gain in $B^r$, see 
\eqref{br}.

Then, modulo acceptable errors (which are zero in the Benjamin-Ono  case),
we have reduced the problem to the case of $R^\sharp$ given by
\[
R^\sharp (\phi, \phi,\phi) = 2 B^\sharp(\phi,\phi \phi_x) - 2 \phi B^\sharp(\phi ,\phi_x) - 2i \tilde C^a(\phi, B^\sharp(\phi,\phi)).
\]

As written, the symbol of $R^\sharp$ has the form
\[
-4i r^\sharp(\xi,\eta,\zeta)
=  \ 2 \Tau(\xi) \Tau(\eta+\zeta) \zeta - 2 \Tau(\eta) \Tau(\zeta) \zeta 
- 2 \tilde C^a(\xi, \eta+\zeta) ( 3 -  \Tau(\eta) \Tau(\zeta))
\]
which we expand as
\[\begin{aligned}
-4i r^\sharp(\xi,\eta,\zeta)
  =& \ 2 \Tau(\xi) \Tau(\eta+\zeta) \zeta - 2 \Tau(\eta) \Tau(\zeta) \zeta - 
\\ &  - \frac12\left[ \xi \Tau (\eta+\zeta)(\Tau(\xi+\eta+\zeta) - \Tau(\xi)) - (\eta+\zeta) 
\Tau (\xi) (\Tau(\xi+\eta+\zeta) - \Tau (\eta+\zeta))\right]
\\ & \ \ \ \ \ \ \ 
 ( 3 -  \Tau(\eta) \Tau(\zeta)).
\end{aligned}
\]
Since the three arguments of $r$ are identical, one should symmetrize the above
expression. However, as an alternative,
we can put all derivatives in $\zeta$ and symmetrize the remaining expression in $\xi$ and $\eta$,
\[
-4i r^\sharp(\xi,\eta,\zeta) = \zeta \Omega(\xi,\eta,\zeta),
\]
where we compute
\[
\begin{aligned}
\Omega(\xi,\eta,\zeta) = & \  \Tau(\xi) \Tau(\eta+\zeta) + \Tau(\eta) 
\Tau(\xi+\zeta) -  \Tau(\eta) \Tau(\zeta) - \Tau(\xi) \Tau(\zeta) 
\\ &
- \frac12 \Tau (\eta+\xi)(\Tau(\xi+\eta+\zeta) - \Tau(\zeta)) 
 ( 3 -  \Tau(\eta) \Tau(\xi))
 \\ & 
  + \frac12 \Tau (\xi) (\Tau(\xi+\eta+\zeta) - \Tau (\eta+\zeta))
 ( 3 -  \Tau(\eta) \Tau(\zeta))
\\ &
+ \frac12 \Tau (\eta) (\Tau(\xi+\eta+\zeta) - \Tau (\xi+\zeta))
 ( 3 -  \Tau(\xi) \Tau(\zeta)).
 \end{aligned}
 \]
To simplify this expression, we use the identity 
\[
\Tau(\xi+\eta) \Tau(\xi) \Tau(\eta) 
=    \Tau(\xi+\eta) - \Tau(\xi) - \Tau(\eta)
\]
as well as the related auxiliary function
\[
G(\xi,\eta) = \Tau(\xi+\eta) \Tau(\xi) +
\Tau(\xi+\eta) \Tau(\eta) -\Tau(\xi) \Tau(\eta) +1 ,
\]
which is easily seen to be a Schwartz function, 
with exponential decay at infinity.

We can use the function $G$ in the last two terms in $\Omega$ to rewrite this as 
\[
\begin{aligned}
\Omega(\xi,\eta,\zeta) \approx & \  \Tau(\xi) \Tau(\eta+\zeta) + \Tau(\eta) 
\Tau(\xi+\zeta) -  \Tau(\eta) \Tau(\zeta) - \Tau(\xi) \Tau(\zeta) 
\\ &
- \frac12 (\Tau(\xi+\eta+\zeta) - \Tau(\zeta)) \Tau (\eta+\xi)
 ( 3 -  \Tau(\eta) \Tau(\xi))
 \\ & 
  - \frac12 (\Tau(\xi+\eta+\zeta) \Tau (\eta+\zeta)+1)
 ( 3 -  \Tau(\eta) \Tau(\zeta))
\\ &
- \frac12  (\Tau(\xi+\eta+\zeta)\Tau (\xi+\zeta) + 1)
 ( 3 -  \Tau(\xi) \Tau(\zeta)).
 \end{aligned}
 \]
where the error terms in $r^\sharp$,
\[
-\frac12 \zeta( G(\xi,\eta+\zeta) ( 3 -  \Tau(\eta) \Tau(\zeta)) + G(\eta,\xi+\zeta) ( 3 -  \Tau(\xi) \Tau(\zeta))
\]
are not directly decaying but do decay after symmetrization. Discarding these exponentially decaying errors, after algebraic simplifications we arrive at
\[
\begin{aligned}
\Omega(\xi,\eta,\zeta) \approx & \  \Tau(\xi) \Tau(\eta+\zeta) + \Tau(\eta) 
\Tau(\xi+\zeta) + \frac12 \Tau(\eta) \Tau(\zeta) + \frac12 \Tau(\xi) \Tau(\zeta) - 3 
\\ &
- \frac12 (\Tau(\xi+\eta+\zeta) - \Tau(\zeta)) \Tau (\eta+\xi)
 ( 3 - \Tau(\eta) \Tau(\xi))
 \\ & 
  - \frac12 \Tau(\xi+\eta+\zeta) \Tau (\eta+\zeta)
 ( 3 -  \Tau(\eta) \Tau(\zeta))
\\ &
- \frac12  \Tau(\xi+\eta+\zeta)\Tau (\xi+\zeta) 
 ( 3 -  \Tau(\xi) \Tau(\zeta)).
 \end{aligned}
 \]
Now we apply the exact $\Tau$ identity in the last three lines, to get 
\[
\begin{aligned}
\Omega(\xi,\eta,\zeta) \approx & \  \Tau(\xi) \Tau(\eta+\zeta) + \Tau(\eta) 
\Tau(\xi+\zeta) -  \frac12 \Tau(\eta) \Tau(\zeta) - \frac12 \Tau(\xi) \Tau(\zeta) + 3 
\\ &
- \frac12 (\Tau(\xi+\eta+\zeta) - \Tau(\zeta))(2 \Tau (\eta+\xi)
  +  \Tau(\eta) + \Tau(\xi))
 \\ & 
  - \frac12 \Tau(\xi+\eta+\zeta) (2 \Tau (\eta+\zeta)
 +  \Tau(\eta) + \Tau(\zeta))
\\ &
- \frac12  \Tau(\xi+\eta+\zeta)(2 \Tau (\xi+\zeta) 
 + \Tau(\xi)+ \Tau(\zeta)).
 \end{aligned}
 \]
We reorganize this as 
\[
\begin{aligned}
\Omega(\xi,\eta,\zeta) \approx & \  \Tau(\xi) \Tau(\eta+\zeta) + \Tau(\eta) 
\Tau(\xi+\zeta) + \Tau(\zeta)\Tau(\xi+\eta)
 + 3 
\\ &
- \Tau(\xi+\eta+\zeta) (\Tau (\eta+\xi)+\Tau(\zeta))
 \\ & 
  -  \Tau(\xi+\eta+\zeta) ( \Tau (\eta+\zeta) +\Tau(\xi))
\\ &
- \Tau(\xi+\eta+\zeta)(\Tau (\xi+\zeta) +\Tau(\eta)) 
\\
= & \  -G(\eta+\xi,\zeta) - G(\eta+\zeta,\xi) - G(\xi+\zeta,\eta).
 \end{aligned}
 \]
The above combination of $G$ functions is symmetric,
therefore we can also symmetrize the $\zeta$ prefactor
of $\Omega$ to $\xi+\eta+\zeta$, which yields the desired exponential decay.

\end{proof}

We now use the symbol $B$ constructed in the previous proposition
in order to prove the energy bound for the corresponding $L^{NL} \phi$.

\begin{proposition}\label{p:LNL-ee}
 Let $L^{NL} \phi$ be defined as in \eqref{def:LNL}, with $B$ given by Proposition~\ref{p:B}. Let $u$ be a solution to \eqref{ilw} in a time interval $[0,T]$ with initial data satisfying \eqref{small-data}.
 Assume in addition that $u$ satisfies the bootstrap assumption 
 \eqref{boot} and that $T \ll_C \epsilon^{-2}$. Then $u$ 
 satisfies the vector field bound \eqref{LNL-ee} in $[0,T]$.
\end{proposition}

\begin{proof}
We recall that $v$ solves 
the equation
\[
P v = C(\phi,v) + t R(\phi,\phi,\phi) + D(\phi,\phi),
\]
where
\begin{equation}\label{R-defn}
R(\phi,\phi,\phi) = 2 B(\phi,\phi\phi_x) - C(\phi, B(\phi,\phi)).
\end{equation}
Then we have
\[
\begin{aligned}
\frac{d}{dt} \frac12 \||\Tau|^\frac12 v\|_{L^2}^2 = & =
\int |\Tau| v \cdot C(\phi,v) dx +  t \int  |\Tau|^\frac12 v \cdot |\Tau|^\frac12 R(\phi,\phi,\phi)\, dx  + \int |\Tau|^\frac12 v \cdot |\Tau|^\frac12 D(\phi,\phi) \, dx  
\\
:= & \ I_1 + I_2+I_3
\end{aligned}
\]
We separately estimate each of the three terms. For $I_1$ we symmetrize,
\[
I_1 = \frac12 \int  |\Tau|^\frac12 v \cdot  \tilde C_{sym}(\phi,|\Tau|^\frac12 v) \,  dx, \qquad \tilde C_{sym} =  |\Tau|^\frac12 C_\phi |\Tau|^{-\frac12} + |\Tau|^{-\frac12} C_\phi^* |\Tau|^{\frac12}  
\]
To estimate this term we will prove the following bound:
\begin{equation}\label{est-I1}
\| \tilde C_{sym}(u,w)\|_{L^2} \lesssim \|\Tau u\|_{L^\infty} \| w\|_{L^2}    
\end{equation}
For the two remaining terms, we will establish the following bounds:
\begin{equation}\label{est-I2}
\| |\Tau|^\frac12 B(u,u,u)\|_{L^2} \lesssim \|\Tau u\|_{L^\infty} \| u\|_{L^\infty} \| u\|_{L^2}  ,  
\end{equation}
\begin{equation}\label{est-I3}
\| |\Tau|^\frac12 D(u,u)\|_{L^2} \lesssim \||\Tau|^\frac12  u\|_{L^\infty} \| u\|_{L^\infty}   .  
\end{equation}
Assuming these three bounds, we first conclude the proof of the proposition.
For this we use our bootstrap bounds to estimate
\[
\| |\Tau|^\frac12 B(\phi,\phi,\phi)\|_{L^2} \lesssim \epsilon^3 t^{-1}, \qquad \| |\Tau|^\frac12 D(\phi,\phi)\|_{L^2} ]\lesssim \epsilon^2 t^{-\frac12}.
\]
Then our energy relation reads 
\[
\frac{d}{dt} \frac12 \||\Tau|^\frac12 v\|_{L^2}^2  \lesssim \epsilon t^{-\frac12} \||\Tau|^\frac12 v\|_{L^2}^2 
+ ( \epsilon^3 + \epsilon^2 t^{-\frac12} )\||\Tau|^\frac12 v\|_{L^2} 
\]
At the initial time we start with 
\[
\||\Tau|^\frac12 v(0)\|_{L^2} \lesssim \epsilon .
\]
Then a straightforward Gronwall inequality leads to 
\begin{equation}
\||\Tau|^\frac12 v(t)\|_{L^2} \lesssim \epsilon, \qquad t \lesssim \epsilon^{-2},
\end{equation}
as needed. It remains to prove the bounds \eqref{est-I1},
\eqref{est-I2} and \eqref{est-I3}.

\bigskip

\emph{Proof of \eqref{est-I1}:}
The symbol of $\tilde C_{sym}$ is 
\[
\tilde c_{sym}(\xi,\eta) = |\Tau(\zeta)|^\frac12 c(\xi,\eta) |\Tau(\eta)|^{-\frac12} + |\Tau(\zeta)|^{-\frac12} c(\xi,\zeta) |\Tau(\eta)|^{\frac12}, \qquad \xi+\eta+\zeta = 0.
\]
We recall that 
\[
c(\xi,\eta) = i \eta - 2i \tilde c_a(\xi,\eta),
\]
where
\[
\tilde c^a (\xi,\eta) = - \frac14  \xi \Tau (\eta)(\Tau(\xi+\eta) - \Tau(\xi)) +\frac14 \eta \Tau (\xi) (\Tau(\xi+\eta) - \Tau (\eta)).
\]

Roughly speaking, we will show that the symbol $\tilde c_{sym}$ has size 
\[
|\tilde c_{sym}| \lesssim |\Tau(\xi)| e^{-c \min\{|\xi|,|\eta|,|\zeta|\}},
\]
and sufficient regularity for the bound \eqref{est-I1}.
For clarity we separate the analysis into several regions:

\medskip

(i) The region $|\xi|+|\eta|+|\zeta| \lesssim 1$. 
Here we first examine the contribution of $\tilde c^a$, where there are no cancellations. Then modulo smooth factors
this contribution has the form $\xi |\eta|^\frac12 |\zeta|^\frac12$, and the bound \eqref{est-I1} is straightforward.
We next consider the contribution of the $\eta$ term in $C$. Harmlessly replacing the $\Tau$ functions with absolute values,
we get a symbol of the form
\[
|\eta|^\frac12 |\zeta|^\frac12 (\sgn \eta + \sgn \zeta).
\]
This has a commutator structure, precisely it corresponds to the operator 
\[
|D|^\frac12 [H,\phi] |D|^\frac12
\]
whose norm in $L^2$ is bounded by $\| \partial \phi\|_{L^\infty}$. To see this we interpolate between the bound for $[H, \phi]\partial $ from  $L^2$ to $L^2$ and its dual. This goes back to Calder\'on's work on commutator estimates \cite{MR177312}. Alternatively one can prove this bound directly  using the standard Littlewood-Paley trichotomy.

\medskip 

(iii) The region $|\eta| \ll 1 \ll |\zeta|$ (and by duality, the region $|\zeta| \ll 1 \ll |\eta|$).
Here we also have $|\xi| \gg 1$. Here there is no cancellation between the two terms in $\tilde c_{sym}$
so we consider only the first one, discarding the harmless $|\Tau(\zeta)|^\frac12$ factor.
The $\eta$ term in $c$ yields $\eta |\Tau(\eta)|^\frac12$ which gives a bounded multiplier at low frequency.
In the first $\tilde c_a$ term we have exponential decay as $\xi \to \infty$ from the last factor, which defeats the growth
in the $\xi$ factor and the second $\tilde c_a$ is directly bounded.

\medskip

(ii) The region $ 1 \ll |\eta|, |\zeta|$. Here we can discard the $|\tau(\eta)|^{-\frac12}$ and the $|\Tau(\xi)|^{-\frac12}$ 
factors, which leave us with the simplified form
\[
\begin{aligned}
\tilde c_{sym}(\xi,\eta) = & \ |\Tau(\zeta)| c(\xi,\eta)  +  c(\xi,\zeta) |\Tau(\eta)|
\\
= & \ \frac12 (|\Tau(\zeta)| + |\Tau(\eta)|)(c(\xi,\eta)  +  c(\xi,\zeta)) + \frac12 (|\Tau(\zeta)| - |\Tau(\eta)|)(c(\xi,\eta)  -  c(\xi,\zeta))
\end{aligned}
\]
In the first term on the right we get exponential decay from Lemma~\ref{l:C-sym}.
In the second term we get both a $\xi$ factor and exponential decay at infinity from the difference $|\Tau(\zeta)| - |\Tau(\eta)|$
provided that $|\xi| \lesssim |\eta|, |\zeta|$. 

So it remains to consider the case when $1 \ll |\eta| \ll |\xi|, |\zeta|$
(and the symmetric case with $\eta$ and $\zeta$ interchanged). Here by Lemma~\ref{l:C-sym} we can replace the difference 
$(\xi,\eta)  -  c(\xi,\zeta)$ by $c(\xi,\eta)$. The difference $|\Tau(\zeta)| - |\Tau(\eta)|$ still decays exponentially in $\eta$ 
and thus will control the $\eta$ terms in $c(\xi,\eta)$. We are left with the $\xi$ term in $\tilde c_a(\xi,\eta)$. But this has the factor 
$\Tau(\xi+\eta) - \Tau(\eta)$ which in this case has exponential decay in $\xi$.

\bigskip

\emph{Proof of \eqref{est-I2}:}
We simply discard the $|\Tau|^\frac12$, which has a bounded symbol.
Then we use the symbol bound \eqref{b-size}, which shows that $R$ is essentially of the form $ \phi \cdot \phi \cdot \Tau \phi$. Precisely, by Lemma~\ref{l:R} we can represent  $R(\phi,\phi,\phi)$ in the form
\[
R(\phi,\phi,\phi) = L( \Tau \phi, \phi, \phi) ,
\]
where the $L$ form has an integrable, rapidly decreasing kernel.
Now we use our bootstrap assumption \eqref{boot}
for the first two entries,
\[
|\phi| \leq C \epsilon \omega_0, \qquad 
|\Tau \phi| \leq C \epsilon \omega_1.
\]
Here $\omega_0$ and $\omega_1$ are positive nonconstant functions, but they are slowly varying on the unit scale spatial scale. Since 
\[
\omega_0 \omega_1 \lesssim t^{-1}
\]
it follows that we can estimate
\[
\|  R(\phi,\phi,\phi)\|_{L^2}  \lesssim 
C^2 \epsilon^2 t^{-1} \| \phi\|_{L^2}\lesssim C^2 \epsilon^3 t^{-1},
\]
which suffices exactly up to cubic time.

\bigskip

\emph{Proof of \eqref{est-I3}:}
The contribution of $D$ is also easily estimated, using 
the symbol bound \eqref{d-size}. Precisely, we claim that we have the bound
\begin{equation}\label{D-est}
\||\Tau|^\frac12 D(\phi,\phi)\|_{L^2} \lesssim  
C \epsilon^2 t^{-\frac12}.
\end{equation}
To prove this we consider a Littlewood-Paley decomposition of each of the two factors. If one is at frequency $>0$ then we directly estimate
\[
\Vert \Tau^{\frac{1}{2}}D(\phi_{>0} , \phi)\Vert_{L^2} \lesssim \Vert \phi_{>0 }\Vert_{L^{\infty}} \Vert \phi\Vert_{L^2} \lesssim 
C \epsilon^2 t^{-\frac12}.
\]
Now it remains to investigate the case when both frequencies are $<0$, where  we use the Besov bound \eqref{besov-u} and the bootstrap bound \eqref{boot-d}
\[
\begin{aligned}
\Vert \Tau^{\frac{1}{2}}D(\phi_{<0},\phi_{<0})\Vert_{L^2}& \lesssim  \sum_{j \leq k < 0} 
\Vert \Tau^{\frac{1}{2}}D(\phi_j,\phi_k)\Vert_{L^2}
\\ 
& \lesssim  \sum_{j \leq k < 0} 
 2^{\frac{k}{2}}
\|\phi_j\|_{L^2} \|\phi_k\Vert_{L^\infty}
\\
&\lesssim  \sum_{j \leq k < 0} \epsilon^2C 2^{\frac{j}{2}}t^{-\frac{1}{2}}\\
&\lesssim  C\epsilon^2 t^{-\frac{t}{2}}.
\end{aligned}
\]

\end{proof}

To switch from the energy bound \eqref{LNL-ee} for $L^{NL} \phi$  to 
the energy bound for $L \phi$ it suffices to estimate the $tB(\phi,\phi)$ correction perturbatively:

\begin{lemma}
We have 
\begin{equation}
 \||\Tau|^\frac12 B(\phi,\phi) \|_{L^2} \lesssim \|  |\Tau|^\frac12 \phi\|_{L^\infty} \|\phi\|_{L^2}  .
\end{equation}   
\end{lemma}

\begin{proof}

 This is similar to \eqref{est-I3}.  
\end{proof}

Once we have this, we obtain the desired bound 
\begin{equation}
\| |\Tau|^\frac12 L \phi\|_{L^2} \lesssim \epsilon, \qquad t \ll \epsilon^{-2}.   
\end{equation}

\vspace{1in}

\bibliography{quasi}

\bibliographystyle{plain}

\end{document}